\documentclass[10pt,leqno]{article}
\usepackage{amssymb,amsbsy,amsmath,amsfonts,amssymb,amscd, mathrsfs}

\usepackage[english]{babel}
\usepackage[T1]{fontenc}
\usepackage{indentfirst}

\makeatletter
\@addtoreset{equation}{section}
\makeatother

\newtheorem{statement}{}[section]
\newtheorem{theoreme}[statement]{Theorem}
\newtheorem{lemme}[statement]{Lemma}

\newtheorem{proposition}[statement]{Proposition}

\newtheorem{corollaire}[statement]{Corollary}

\newcommand\C{\mathbb C}

\newcommand\R{\mathbb R}
\newcommand\T{\mathbb T}
\newcommand\D{\mathbb D}

\newcommand\e{{\rm e}}

\newcommand\ind{{\rm 1\kern-.30em I}}
\newcommand\qed{\hfill $\square$}
\renewcommand \Re{{\mathfrak R}{\rm e}\,}
\renewcommand \Im{{\mathfrak I}{\rm m}\,}

\title{\bf   On approximation numbers of composition operators }

\author{\it  Daniel Li,  Herv\'e Queff\'elec, Luis Rodr{\'\i}guez-Piazza} 

\date{\footnotesize \today}

\begin{document}

\maketitle

\noindent{\bf Abstract.} \emph{We show that the approximation numbers of a compact composition operator
on  the weighted Bergman spaces $\mathfrak{B}_\alpha$ of the unit disk can tend to 0 arbitrarily slowly, but that they never tend quickly to 0: they grow at 
least exponentially, and this speed of convergence is only obtained for symbols which do not approach the unit circle. We also give an upper bounds and 
explicit an example.} 
\medskip

\medskip



\section{Introduction} \label{section 1}

Let $\D$ be the open unit disk of the complex plane, equipped with its normalized area measure $dA (z) = \frac{dxdy}{\pi} \,$.
For $\alpha > - 1$, let $\mathfrak{B}_\alpha$ be the weighted Bergman space of analytic functions $f (z) =\sum_{n=0}^\infty a_n z^n$ on $\D$ such that
\begin{displaymath}
\Vert f \Vert_{\alpha}^2 = (\alpha + 1)\int_{\D} \vert f (z)\vert^2 (1 - \vert z \vert^2)^{\alpha} \,dA (z) 
= \sum_{n=0}^\infty \frac{n! \Gamma (2 + \alpha)}{\Gamma (n + 2 + \alpha)} \vert  a_n \vert^2 < \infty. 
\end{displaymath}
The limiting case, as $\alpha \mathop{\longrightarrow}\limits^{>} - 1$, of those spaces is the usual Hardy space $H^2$ (indeed, if $f$ is a polynomial, we 
have 
$\lim_{\alpha \mathop{\longrightarrow}\limits^{\scriptscriptstyle >} - 1} \Vert f\Vert_{\alpha}^2 =\sum_{n=0}^\infty \vert a_n\vert^2 
= \Vert f\Vert_{H^2}^2$), 
which we shall treat as $\mathfrak{B}_{-1}$. Note that $\Vert f \Vert_{\alpha}^2\approx \sum_{n=0}^\infty \frac{\vert a_n\vert^2}{(n+1)^{\alpha+1}}$  
and that  
\begin{displaymath}
dA_{\alpha}(z) = (\alpha + 1) (1 - \vert z\vert^2)^\alpha \, dA (z) 
\end{displaymath}
is a probability measure on $\D$. 
\smallskip

Bergman spaces  (\cite{Zhu-livre} page~75, page~78) are Hilbert spaces of analytic functions on $\D$ with reproducing kernel $K_a \in \mathfrak{B}_\alpha$, 
given by $ K_{a}(z)=(\frac{1}{1 -\overline{a}z})^{\alpha+2}$, namely, for every $a\in\D$:
\begin{equation}\label{Bergman} 
f (a) = \langle f, K_a\rangle\,, \  \forall f \in \mathfrak{B}_\alpha ; \quad \text{and} \quad 
\Vert K_a\Vert^2 = K_{a}(a) = \smash{\Big(\frac{1}{1-\vert a\vert^2}\Big)^{\alpha + 2} } .
\end{equation}
An important common feature of those spaces is that the  multipliers of $\mathfrak{B}_\alpha$ can be (isometrically) identified with the space $H^\infty$ 
of bounded analytic functions on $\D$, that is: 
\begin{equation}\label{Multiplier} 
\forall g\in H^\infty,\qquad  \Vert g\Vert_\infty =\sup_{f\in \mathfrak{B}_\alpha, \Vert f \Vert_{\alpha}\leq 1}\Vert fg\Vert_\alpha.
\end{equation}
Indeed,  $\Vert f g \Vert_\alpha \leq \Vert g \Vert_\infty \Vert f \Vert_\alpha$ is obvious, and if 
$\Vert f g \Vert_\alpha \leq C \Vert f \Vert_\alpha$ for\ all $f\in \mathfrak{B}_\alpha$, testing this inequality successively on $f = 1, g, \ldots, g^n,\ldots$  
easily gives $g\in H^\infty$ and $\Vert g \Vert_\infty\leq C$.
\medskip

Let now $\varphi$ be a \emph{non-constant} analytic self-map (a so-called \emph{Schur function}) of  $\D$ and let 
$C_\varphi \colon \mathfrak{B}_\alpha\to {\mathscr H} (\D)$ the associated  composition operator:
\begin{displaymath}
C_\varphi (f) = f \circ \varphi.
\end{displaymath}
It is well-known (\cite{COMC-livre} page~30) that such an operator is always bounded from $\mathfrak{B}_\alpha$ into itself, and we are interested in its 
approximation numbers.
\par\medskip

Also recall that the \emph{approximation} (or \emph{singular}) \emph{numbers} $a_{n} (T)$ of an operator $T \in \mathcal{L}(H_1, H_2)$, between two 
Hilbert spaces $H_1$ and $H_2$,  are defined, for $n = 1,2,\ldots \,$, by: 
\begin{displaymath}
a_{n}(T) = \inf \{\Vert T - R\Vert \,;\  {\rm rank}\, (R) < n \} .  
\end{displaymath}
We have:
\begin{displaymath}
a_{n}(T) = c_{n}(T) = d_{n}(T) \,,     
\end{displaymath}
where the numbers $c_{n}$ (resp. $d_n$) are the \emph{Gelfand} (resp. \emph{Kolmogorov}) \emph{numbers} of $T$ (\cite{CA-ST-livre}, page~59 and 
page~51 respectively). \par\smallskip

In the sequel we shall need the following quantity:
\begin{equation} \label{beta}
\beta (T) = \liminf_{n\to\infty} \big[ a_{n}(T) \big]^{1/ n} .
\end{equation}
Those approximation numbers form a non-increasing sequence such that 
\begin{displaymath}
a_{1} (T) = \Vert T \Vert, \qquad a_{n}(T) = a_{n} (T^*) = \sqrt{a_{n} (T^*T)} 
\end{displaymath}
and verify the so-called  ``ideal''   and ``subadditivity'' properties (\cite{K\"onig-livre} page~57 and page~68): 
\begin{equation}\label{ideal}
a_{n} (ATB)\leq\Vert A\Vert \, a_{n} (T) \, \Vert B\Vert \,; \quad   a_{n+m-1}(S+T)\leq a_{n}(S)+a_{m}(T) .
\end{equation}       
Moreover,  the sequence $(a_{n}(T))$ tends to $0$ iff $T$ is compact. If $(a_{n} (T)) \in\ell_{p}$, we say that $T$ belongs to  the Schatten class $S_p$ of index 
$p$, $0 < p <\infty$. Taking for $T$ a compact diagonal operator, we see that this sequence  is  non-increasing with limit $0$, but otherwise arbitrary.
But if we restrict ourselves to a \emph{specified class} of operators, the answer is far from being so simple, although in some cases the situation is completely 
elucidated. For example, for the class of  Hankel operators on $H^2$ (those operators $\mathcal{H}_{\phi}$ whose matrix $(a_{i,j})$ on the canonical basis 
of $H^2$ is of the form $a_{i,j} =\widehat\phi (i+j)$ for some function $\phi \in L^\infty$), it is known that $\mathcal{H}_{\phi}$ is compact if and only if 
the conjugate $\bar \phi$ of the symbol $\phi$ belongs to $ H^\infty + {\mathcal C}$, where ${\mathcal C}$ denotes the space of continuous, $2\pi$-periodic 
functions  (Hartman's theorem, \cite{Nikolski-livre 1} page~214). For those Hankel operators,  the following theorem, due to  A.~V.~Megretskii, V.~V.~Peller, 
and S.~R.~Treil  (\cite{MPT} and \cite{Peller-livre}, Theorem~0.1, page~490), shows that the approximation numbers are absolutely  arbitrary, under the 
following form.
\begin{theoreme}[Megretskii-Peller-Treil] \label{Hankel} 
Let  $(\varepsilon_n)_{n\geq 1}$ be a non-increasing sequence of positive numbers. Then there exists a  Hankel operator $\mathcal{H}_{\phi}$ satisfying:
\begin{displaymath}
\qquad a_{n} (\mathcal{H}_{\phi}) = \varepsilon_n, \qquad  \forall n \geq 1.
\end{displaymath}
\end{theoreme}
Indeed, if we take a positive self-adjoint operator $A$ whose eigenvalues $s_n$ coincide with the $\varepsilon_n$'s and whose kernel is infinite-dimensional, 
it is easily checked that this operator $A$ verifies the three necessary and sufficient conditions of Theorem~0.1, page~490 in \cite{Peller-livre} and is therefore 
unitarily equivalent to a Hankel operator $\mathcal{H}_{\phi}$ which will verify, in view of \eqref{ideal}: 
\begin{displaymath}
\qquad \qquad a_{n} (\mathcal{H}_{\phi}) = a_{n} (A) = \varepsilon_n, \qquad n = 1, 2, \ldots 
\end{displaymath}
In particular, if $\varepsilon_n\to 0$, the above Hankel operator will be compact, and in no Schatten class if $\varepsilon_n = 1/\log(n+1)$ for example. 
We also refer to \cite{KP} for the following slightly weaker form due to S.~V.~Khrusc\"ev and V.~Peller, but  with a more elementary proof based on 
interpolation sequences in the Carleson sense: for any $\delta > 0$, there exists a Hankel operator $\mathcal{H}_{\phi}$ such that 
\begin{displaymath}
\qquad \qquad \quad \frac{1}{1+\delta} \,\varepsilon_n \leq a_{n} (\mathcal{H}_{\phi}) \leq (1+\delta)\,\varepsilon_n, \qquad  n=1, 2,\ldots
\end{displaymath}

Now, the aim of this work is to prove  analogous theorems  for the class of composition operators (whose compactness was characterized in \cite{McCluer} 
and \cite{Shap}). But if  we are able  to obtain the Khrusc\"ev-Peller analogue for the  \emph{lower bounds}, we  will only obtain subexponential estimates 
for the upper bounds, a fact  which  is explained by our second result: the speed of convergence to $0$ of the approximation numbers of a composition operator 
cannot be greater than geometric (and \emph{is} geometric for symbols $\varphi$ verifying $\Vert \varphi\Vert_\infty < 1$). Our first result involves a constant 
$< 1$  and is not as precise as the result of Megretskii-Peller-Treil or even that of Khrusc\"ev-Peller; this  is apparently due to the  non-linearity of the 
dependence with respect to the symbol for the class of composition operators, contrary to the case of the  Hankel class.This latter lower bound  improves several 
previously known results on ``non-Schattenness'' of those operators (see Corollary~\ref{Coco} below) and also answers in the positive to a  question which was 
first asked to us by C.~Le~Merdy (\cite{Le Merdy}) in the OT Conference 2008 of Timisoara, concerning the bad rate of approximation of compact composition 
operators. Those  theorems are, to our knowledge, the first individual results on approximation numbers $a_{n}$ of composition operators (in the work of 
Parfenov \cite{Parf},  some good estimates are given for the approximation numbers of the Carleson embedding operator in the case of the space 
$H^2 =\mathfrak{B}_{-1}$,  but they remain fairly implicit, and are not connected with composition operators), whereas all previous results where in terms of 
symmetric norms of the sequence $(a_n)$, not on the behaviour of each $a_n$. 
\medskip\goodbreak

Before describing our results, let us recall two definitions.  For every $\xi$ with $\vert\xi\vert=1$ and $0<h<1$, the Carleson window $W(\xi,h)$ centered at 
$\xi$ and of size $h$ is the set
\begin{displaymath}
W (\xi,h) = \{z\in \overline{\D} \,;\ \vert z\vert\geq 1-h \text{ and } \vert \arg (z\overline{\xi}) \vert\leq \pi h\}.
\end{displaymath}
Let  $\mu$ be a positive, finite, measure on $\overline{\D}$; the associated  maximal function $\rho_\mu$ is defined by: 
\begin{equation}\label{Maximal}
\rho_{\mu} (h) = \sup_{\vert \xi\vert = 1} \mu \big(W(\xi,h)\big).
\end{equation}
 The measure $\mu$ is called a \emph{Carleson measure} for the Bergman space $\mathfrak{B}_\alpha$, or an \emph{$(\alpha+2)$-Carleson measure} 
(including the case $\mathfrak{B}_{-1} = H^2$),  if $\rho_{\mu} (h) = O\, (h^{2+\alpha})$ as $h\to 0$. For any Schur function $\varphi$, we shall denote by 
$m_\varphi$ the image $\varphi^*(m)$ of the Haar measure $m$ of the unit circle under the radial limits function  
$\varphi^*(u) = \lim_{r\to 1^-} \varphi (ru)$ of $\varphi$, $\vert u \vert = 1$, and by $A_{\varphi,\alpha+2}$ the image of the probability measure 
$(\alpha+1) (1 - \vert z \vert^2)^\alpha dA(z)$ under $\varphi$. The corresponding maximal function will be denoted by $\rho_{\varphi, \alpha+2}$. 
This notation is justified by the fact that  $m_\varphi \mathop{=}\limits^{def} A_{\varphi,1}$ is a $1$-Carleson measure and  $A_{\varphi,\alpha}$ an 
$(\alpha+2)$-Carleson measure for $\alpha> -1$, in view of the famous  Carleson embedding theorem which, expressed under a quantitative and generalized 
form, states the following, implicit as concerns $\Vert j \Vert$ and with different notations, but fully proved in \cite{Steg}, Theorem~1.2, for the case 
$\alpha>-1$ (see \cite{Nikolski-livre 1}, page~153).
\begin{theoreme}[Carleson's theorem] 
For any $(\alpha+2)$-Carleson measure $\mu$, the canonical inclusion mapping $j \colon\mathfrak{B}_\alpha\to L^{2}(\mu)$ is defined and  continuous, 
and its norm satisfies 
\begin{equation}\label{A priori}
C^{-1}\sup_{0 < h < 1} \sqrt{\frac{\rho_{\mu} (h)}{h^{2+\alpha}}} \leq \Vert j \Vert 
\leq C \sup_{0 < h < 1} \sqrt{\frac{\rho_{\mu} (h)}{h^{2+\alpha}}} \cdot
\end{equation}
\end{theoreme}
\bigskip

The paper is organized as follows. Section~\ref{section 1} is this introduction. In Section~\ref{section 2}, we prove some preliminary lemmas. Our first 
theorems concern lower bounds. In Section~\ref{section 3}, we prove (Theorem~\ref{Secondary}) that the convergence of the approximation numbers 
$a_{n} (C_\varphi)$ of a composition operator $C_\varphi \colon {\mathfrak B}_\alpha \to {\mathfrak B}_\alpha$ cannot exceed an exponential speed: 
for some $r \in (0,1)$ and some constant $c > 0$, one has $a_{n} (C_\varphi) \geq c\, r^n$. More precisely, with the notations \eqref{beta} and 
\eqref{crochet}, one has $\beta (C_\varphi) \geq [\varphi]^2$. Moreover, this speed of convergence is only attained if the values of $\varphi$ do not approach 
the boundary of the unit disk: $\| \varphi\|_\infty < 1$ (Theorem~\ref{Lastresult}). On the other hand, the speed of convergence to $0$ of $a_{n} (C_\varphi)$ 
can be arbitrarily slow; this is proved in Section~\ref{section 4}. The proof is mainly an  adaptation of the one in \cite{Carroll-Cowen}, but  is fairly technical at 
some points, and will require several additional  explanations. In Section~\ref{section 5}, we prove an upper estimate (Theorem~\ref{Ternary}), and give three 
applications of this theorem. In the final Section~\ref{section 6}, we test our general results against  the example of lens maps, which are known to generate 
composition operators belonging to all Schatten classes.  


\section{Preliminary lemmas}\label{section 2}
 
In this Section, we shall state several lemmas, which are either already known or quite elementary, but turn out to be necessary for the proofs of  our 
Theorem~\ref{Secondary} and Theorem~\ref{Principal}. 
\smallskip

For the proof of Theorem~\ref{Secondary}, we shall need the Weyl lemma (\cite{CA-ST-livre} Proposition~4.4.2, page~157).
\begin{lemme} [Weyl lemma]\label{Weyl}  
Let $T \colon H\to H$ be a compact operator. Suppose that   $(\lambda_n)_{n\geq 1}$ is the sequence of eigenvalues of $T$ rearranged in non-increasing 
order. Then, we have:
\begin{displaymath}
\prod_{k = 1}^n a_{k} (T) \geq \prod_{k=1}^n \vert \lambda_k\vert.
\end{displaymath}
\end{lemme}
\smallskip

We recall (\cite{Carleson}, \cite{Hoffman-livre} pages~194--195, \cite{Nikolski-livre 2} pages~302--303) that an \emph{interpolation sequence} $(z_n)$ with 
(best) interpolation constant $C$ is a sequence $(z_n)$ (necessarily Blaschke, i.e. $\sum_{n = 1}^\infty (1 - \vert z_n \vert) < \infty$)  in the unit disk such 
that, for any bounded sequence $(w_n)$ of scalars, there exists a bounded analytic function $f$ (i.e. $f\in H^\infty$) such that:
\begin{displaymath}
f (z_n) = w_n \,,\ \forall n \geq 1, \quad \text{and}  \quad \Vert f \Vert_\infty \leq C \sup\nolimits_{n \geq 1} \vert w_n\vert.
\end{displaymath} 
The \emph{Carleson constant} $\delta$ of a Blaschke sequence $(z_n)$ is defined as follows:
\begin{equation}\label{Carleson constant}
\delta_n = \prod_{j\neq n} \rho(z_n, z_j) \,; \qquad \delta = \inf \delta_n = \inf_{n\geq 1} (1-\vert z_n \vert^2) \vert B'(z_n)\vert, 
\end{equation}
where $B$ is the Blaschke product with zeroes $z_n$, $n \geq 1$. The interpolation constant $C$ is related to the Carleson constant $\delta$ by the following 
inequality (\cite{Garnett-livre} page~278), in which $\lambda$ is a positive numerical constant:
\begin{equation}\label{Interpolation constant}
\frac{1}{\delta} \leq C \leq \frac{\lambda}{\delta} \bigg(1+\log \frac{1}{\delta}\bigg) \cdot
\end{equation}
This latter inequality can be viewed as a quantitative form of the Carleson interpolation theorem. Interpolation sequences and reproducing kernels of 
$\mathfrak{B}_\alpha$ are related as follows (\cite{Nikolski-livre 2} pages~302--303).
\begin{lemme}\label{Riesz basis}
Let $(z_n)_{n\geq 1}$ be an $H^\infty$-interpolation sequence of the unit disk, with  interpolation constant $C$. Then, the sequence 
$ (f_n) = (K_{z_n}/ \Vert K_{z_n}\Vert)$ of normalized reproducing kernels at $z_n$ is $C$-equivalent to an orthonormal basis in 
$\mathfrak{B}_\alpha$, namely we have for any finite sequence $(\lambda_n)$ of scalars:
\begin{equation}\label{Equivalent} 
C^{-1}\Big(\sum\nolimits_n \vert \lambda_n \vert^2\Big)^{1/2} \leq \Big\Vert  \sum\nolimits_n \lambda_{n} f_n\Big\Vert_\alpha
\leq C \Big(\sum\nolimits_n \vert \lambda_n\vert^2\Big)^{1/2}.
\end{equation}
\end{lemme}
\smallskip

\noindent
The proof in \cite{Nikolski-livre 2} is only for $H^2$, therefore we indicate a simple proof  valid for Bergman spaces $\mathfrak{B}_\alpha$ as well. Let 
$S = \sum \lambda_{n} K_{z_n}$ be a finite linear combination of the kernels $K_{z_n}$, $\omega = (\omega_n)$ be a sequence of complex signs, 
$S_\omega = \sum \omega_n \lambda_{n} K_{z_n}$ and $g\in H^\infty$ an interpolating function for the sequence 
$(\overline{\omega}_n)$, i.e. $g(z_n) = \overline{\omega}_n$ and $\Vert g\Vert_\infty\leq C$. If $f\in \mathfrak{B}_\alpha$ and 
$\Vert f \Vert_\alpha\leq 1$, we see that:
\begin{displaymath}
\langle S_{\omega}, f \rangle = \sum \omega_{n} \lambda_{n} \overline{f(z_n)} = \sum \lambda_{n}\overline{(fg)(z_n)}  
=\sum \lambda_{n}\langle K_{z_n}, fg \rangle =\langle S, fg\rangle , 
\end{displaymath}
so that using \eqref{Multiplier}:
\begin{displaymath}
\vert \langle S_{\omega}, f\rangle \vert \leq \Vert S\Vert_\alpha \Vert fg\Vert_\alpha 
\leq \Vert S \Vert_\alpha \Vert g \Vert_\infty \Vert f\Vert_\alpha\leq C\Vert S \Vert_\alpha 
\end{displaymath}
and passing to the supremum on $f$, we get $\Vert S_\omega\Vert_\alpha\leq C\Vert S\Vert_\alpha$. Since the coefficients $\lambda_{n}$ are arbitrary, 
this implies that $(f_n)$ is $C$-unconditional, namely:
\begin{displaymath}
C^{-1} \Big\Vert \sum \omega_n \lambda_{n} f_n\Big\Vert_\alpha \leq\Big \Vert \sum  \lambda_{n} f_n\Big\Vert_\alpha 
\leq C \, \Big\Vert \sum \omega_n \lambda_{n} f_n\Big\Vert_\alpha. 
\end{displaymath}
Now, squaring and integrating with respect to random, independent, choices of signs $\omega_{n}$'s, we get \eqref{Equivalent}. \qed
\par\bigskip

We also recall (\cite {Hoffman-livre} pages~203--204) that an increasing sequence  $(r_n)$ of numbers such that $0 < r_n < 1$ and 
$\frac{1 - r_{n+1}}{1 - r_{n}}\leq \rho <1$ (i.e. verifying the so-called \emph{Hayman-Newman condition}) is an interpolation sequence (see also 
\cite{Nikolski-livre 1}). In the following, let $(r_n)$ be such a sequence verifying moreover  the backward induction relation:
\begin{equation}\label{Backward} 
\smash {\varphi(r_{n+1}) = r_n. }
\end{equation}   
Set $f_n = K_{r_n}/ \Vert K_{r_n}\Vert$ and $W = \overline{{\rm span}} (f_n)$.  Let $(e_n)_{n\geq 1}$ be the canonical basis of $\ell^2$, $\varphi$ a 
Schur function and $h \in H^\infty$ a function vanishing at $r_1$. Denote by $M_h  \colon \mathfrak{B}_\alpha \to \mathfrak{B}_\alpha$ the operator of 
multiplication by $h$. Then, we have the following basic lemma, which shows that some compression of $C_{\varphi}^*$ is a  backward shift with controlled 
weights (\cite{Carroll-Cowen}).
\begin{lemme}\label{Compression} 
Let $J \colon \ell^2\to W$ be the isomorphism given by $J(e_n) = f_n$. Then, the operator 
${\mathbf B} = J^{-1} C_{\varphi}^* M_{h}^* J \colon \ell^2\to\ell^2$ is the weighted  backward shift given by:
\begin{equation}\label{Good Weights}
{\mathbf B} (e_{n+1}) = w_n e_n \quad \text{and} \quad {\mathbf B}(e_1) = 0, \quad 
\text{where } \quad w_n = \overline{h (r_{n+1})} \frac{\Vert K_{r_n}\Vert}{\Vert K_{r_{n+1}}\Vert} \,\cdot
\end{equation}
\end{lemme}
\goodbreak

To exploit Lemma~\ref{Compression}, we shall need the following simple fact on approximation numbers of weighted backward shifts. 
\begin{lemme}\label{Simple minoration} 
Let $(e_n)_{n\geq 1}$ be an orthonormal basis of the Hilbert space $H$ and ${\mathbf B} \in \mathcal{L} (H)$ the weighted backward shift defined by 
\begin{displaymath}
\qquad {\mathbf B}(e_1) = 0 \quad \text{and} \quad {\mathbf B} (e_{n + 1}) = w_n e_{n}, \qquad \text{where }  w_n\to 0.  
\end{displaymath}
Assume that  $\vert w_n\vert \geq \varepsilon_{n}$ for all $n \geq 1$,  where $(\varepsilon_n)$ is a non-increasing sequence of positive numbers. Then 
${\mathbf B}$ is compact, and satisfies: 
\begin{equation} \label{mino shift} 
a_{n} ({\mathbf B}) \geq \varepsilon_{n},\quad    \forall n\geq 1.
\end{equation}
\end{lemme}

\noindent{\bf Proof.} The compactness of ${\mathbf B}$ is obvious. Let $R$ be an operator of rank $< n$. Then $\ker R$ is of codimension $< n$, and 
therefore intersects the $n$-dimensional space generated by $e_{2},\ldots, e_{n+1}$ in a vector $x = \sum_{j=1}^{n} x_j e_{j+1}$ of norm one. We then have:
\begin{align*}
\Vert {\mathbf B} - R \Vert^2 
& \geq \Vert {\mathbf B} x - Rx \Vert^2 = \Vert {\mathbf B}x \Vert^2 
=\sum\nolimits_{j=1}^{n} \vert w_{j} \vert^2\vert x_j\vert^2 \\ 
& \geq \sum\nolimits_{j=1}^{n} \varepsilon_{j}^2 \vert x_j \vert^2 
\geq \varepsilon_{n}^2 \sum\nolimits_{j=1}^{n} \vert x_j \vert^2 = \varepsilon_{n}^2.
\end{align*}
This  ends the proof of Lemma~\ref{Simple minoration}. \qed
\bigskip

Now, in view of \eqref{Bergman} and  \eqref{Good Weights}, the weight $w_n$ roughly behaves as 
$\sqrt{\frac{1 - r_{n+1}}{1 - r_n}}$, so we shall need good estimates on that quotient, before defining the sequence $(r_n)$  explicitly.
\medskip

 We first connect this estimate with the hyperbolic distance $d$ in $\D$. We denote (see \cite{Hayman} or \cite{Keen-livre} for the definition) by  
$d(z, w; U)$  the  hyperbolic distance of two points  $z, w$ of a simply connected domain $U$.  It follows from the generalized Schwarz-Pick lemma 
(\cite{Keen-livre} Theorem~7.3.1, page~130) applied to the canonical injection $U\to V$ that the bigger the domain the smaller the hyperbolic distance, namely:
 \begin{equation}\label{Principle}
U\subset V \text{ and } z, w\in U \quad \Longrightarrow \quad d (z, w; V) \leq d (z, w; U).
\end{equation}
 Moreover, as is well-known, 
\begin{displaymath}
\smash {0 \leq r < 1 \quad \Longrightarrow \quad d (0, r; \D) = \frac{1}{ 2} \log \frac{1 + r}{1 - r} \,\cdot}
\end{displaymath}
Recall that the pseudo-hyperbolic and hyperbolic distances $\rho$ and $d$ on $\D$ are defined by:
\begin{displaymath}
\rho (a, b) =\Big\vert \frac{a - b}{1-\overline{a}b} \Big\vert \, \raise 1,5pt \hbox{,} \qquad d(a, b) = \frac{1}{2}\log \frac{1+\rho(a,b)}{1 - \rho(a, b)} 
\raise 1pt \hbox{,} \quad  a, b\in\D, 
\end{displaymath}
In the sequel, we shall omit the symbol $\D$ as far as the open unit disk is concerned. For this unit disk, we have the following simple inequality 
(\cite{Carroll-Cowen}) .
\begin{lemme}\label{Poincare} 
Let $ a, b\in \D$ with $0 < a < b < 1$. Then: 
\begin{equation}\label{Hyperbolic}
\e^{- 2 d(a, b)} \leq \frac{1 - b}{1 - a} \leq 2 \, \e^{- 2 d(a, b)}. 
\end{equation}
\end{lemme}
\goodbreak 

Finally, before proceeding to the construction of our Schur function $\varphi$ in Section 4, it will be useful to note the following simple technical lemma.
\begin{lemme} \label{Logconvex} 
Let $(\varepsilon_n)$ be a non-increasing sequence of positive numbers of limit $0$. Then there exists a decreasing and logarithmically convex sequence 
$(\delta_n)$ of positive numbers, with limit $0$, such that $\delta_n \geq \varepsilon_n$ for all $n\geq 1$.
 \end{lemme}
 
\noindent{\bf Proof.} Provided that we replace $\varepsilon_n$ by $\varepsilon_n + \frac{1}{n}$, we may assume that $(\varepsilon_n)$ is decreasing. Let us 
define our new sequence by the inductive relation: 
\begin{displaymath}
\smash{ \delta_1 = \varepsilon_1\,;\quad  \delta_2 = \varepsilon_2\,; \quad 
\delta_{n+1} = \max\big(\varepsilon_{n+1}, \delta_{n}^2 / \delta_{n-1}\big) .}
\end{displaymath}
 This sequence is log-convex by definition, i.e. $\delta_{n}^2 \leq \delta_{n+1}\delta_{n-1}$. By induction, it is seen to be decreasing. Therefore, it has a 
limit $l \geq 0$. If $\delta_n =\varepsilon_n$ for infinitely many indices, $l=0$. Otherwise, for $n$ large enough, we have the inductive relation 
$\delta_{n+1} = \delta_{n}^2/ \delta_{n - 1}$, which implies that $\delta_n = \exp (\lambda n +\mu)$ for some constants $\lambda, \mu$. Since 
$(\delta_n)$ is decreasing, we must have $\lambda < 0$ and again we get $l = 0$.  \qed
\par\medskip
 
In the sequel, we may and will thus assume, without loss of generality, that $(\varepsilon_n)$ is decreasing and logarithmically convex.  


\section{Lower bounds}\label{section 3} 

We first introduce a notation. If
\begin{displaymath}
\varphi^{\#}(z) =\lim_{w\to z} \frac{\rho (\varphi (w),\varphi(z))}{\rho(w,z)} = 
\frac{\vert \varphi '(z)\vert(1 -\vert z\vert^2)}{1 -\vert \varphi (z)\vert^2} 
\end{displaymath}
is the pseudo-hyperbolic derivative of $\varphi$, we set: 
\begin{equation} \label{crochet} 
[\varphi] = \sup_{z\in\D} \varphi^{\#}(z) =\Vert \varphi^{\#}\Vert_\infty. 
\end{equation}

In our first theorem, we get that the approximation numbers cannot supersede a geometric speed.

\begin{theoreme}\label{Secondary} 
For any Schur function $\varphi$, there exist positive  constants $c > 0$ and $0 < r <1$ such that, for 
$C_\varphi:\mathfrak{B}_\alpha\to\mathfrak{B}_\alpha$, we have: 
\begin{equation}
\qquad \qquad a_{n} (C_{\varphi}) \geq c \,r^n, \qquad n = 1, 2, \ldots
\end{equation}
More precisely, one has $\beta (C_\varphi)  \geq [\varphi]^2$ and hence, for each $\kappa < [\varphi]$, there exists a constant $c_\kappa > 0$ such that:
\begin{equation} \label{Bloch} 
a_{n} (C_\phi)\geq c_\kappa  \kappa^{2n}.
\end{equation}
\end{theoreme}
\goodbreak

For the proof,  we need  the following lemma.
\begin{lemme}\label{Second} 
Let $T \colon H\to H$ be a compact operator. Suppose that   $(\lambda_n)_{n\geq 1}$, the sequence of eigenvalues of $T$ rearranged in non-increasing order, 
satisfies, for some $\delta > 0$ and $r \in (0,1)$:
\begin{displaymath}
\qquad \vert \lambda_n \vert \geq \delta r^n, \quad  n = 1, 2, \ldots 
\end{displaymath}
Then there exists $\delta_1 > 0$ such that  
\begin{displaymath}
\qquad \quad a_{n}(T) \geq \delta_1 r^{2n}, \quad n = 1, 2, \ldots 
\end{displaymath}
In particular $\beta (T) \geq r^2$.
\end{lemme}

\noindent{\bf Proof.} By Weyl's inequality (Lemma~\ref{Weyl}), we have 
\begin{displaymath}
\prod_{k=1}^n a_{k} (T) \geq \prod_{k=1}^n \vert \lambda_k \vert \geq \delta^n r^{n(n+1)/2}. 
\end{displaymath}
Since $a_{k} (T)$ is non-increasing and  $a_{k} (T)\leq \Vert T \Vert$ for every $k$, changing $n$ into $2n$,  we get:
\begin{displaymath}
\Vert T \Vert^n a_{n} (T)^n \geq \prod_{k=1}^{2n} a_{k} (T)\geq \delta^{2n} r^{n(2n+1)} \geq \delta^{2n} r^{2n^2}  
\end{displaymath}
and therefore $ a_{n} (T) \geq  \frac{\delta^2}{\Vert T \Vert} \, r^{2n} =\delta_1 r^{2n}$, as claimed. \qed
\bigskip

By applying this lemma to composition operators, we get the following result, which ends the proof of Theorem~\ref{Secondary}. 

\begin{proposition}\label{Canon} 
For every composition operator $C_\varphi \colon \mathfrak{B}_\alpha\to \mathfrak{B}_\alpha$ of symbol $\varphi \colon\D \to \D$,  
we have $\beta(C_\varphi)\geq [\varphi]^2$.
\end{proposition}
\noindent{\bf Proof.} For every $a\in \D$, let $\Phi_a$ be the (involutive) automorphism of the unit disk defined by 
\begin{displaymath}
\qquad \quad \Phi_{a} (z) = \frac{a - z}{1 - \overline{a}z} \, \raise 1,5pt \hbox{,} \quad z \in \D.
\end{displaymath}
Observe that we have 
\begin{displaymath}
\Phi_{a} (a) = 0, \quad \Phi_{a} (0) = a, \quad \Phi'_{a} (a) = \frac{1}{\vert a \vert^2 - 1} \raise 1,5pt \hbox{,} \quad \Phi'_{a} (0) = \vert a \vert^2 - 1. 
\end{displaymath}
Define now $\psi = \Phi_{\varphi(a)} \circ \varphi \circ \Phi_a$.  We have that $0$ is a fixed point of $\psi$, whose derivative is, by the chain rule:
\begin{equation} \label{Dieze} 
\psi' (0) = \Phi'_{\varphi (a)}(\phi (a)) \varphi' (a) \Phi'_{a} (0) = 
\frac{\varphi' (a) (1 - \vert a \vert^2)}{1 - \vert \varphi (a) \vert^2} \mathop{=}^{def} \varphi^{\#} (a).
\end{equation}
By Schwarz's lemma, we know that $\vert \psi' (0) \vert \leq 1$ and so 
$\frac{\vert \varphi' (a) \vert (1 - \vert a \vert^2)}{1 - \vert \varphi (a) \vert^2} \leq 1$ (Schwarz-Pick's inequality). \par
Let us first assume that the composition operator $C_\varphi$ is compact. Then, so is  $C_\psi$, since we have 
\begin{equation} \label{Equivalent 2}
C_\psi = C_{\Phi_a}\circ C_\varphi \circ C_{\Phi_{\varphi (a)}}. 
\end{equation}
If $\psi' (0) \neq 0$, the sequence of  eigenvalues of $C_\psi$  is $\big( [\psi' (0)]^n\big)_{n\geq 0}$ (\cite{Shap-livre}, page~96; the result given for the space 
$H^2$ holds for $\mathfrak{B}_\alpha\subset H^2$, and would also hold for \emph{any} space of analytic functions in $\D$ on which $C_\psi$ is compact). 
Lemma~\ref{Second} then gives us: 
\begin{displaymath}
\beta (C_\psi) \geq \vert \psi' (0) \vert^2 = [\varphi^{\#} (a) ]^2 \geq 0.
\end{displaymath}
This trivially still holds if  $\psi' (0) = 0$. \par
Now, since $C_{\Phi_a}$ and $ C_{\Phi_{\varphi(a)}}$ are invertible operators, \eqref{Equivalent 2} clearly implies that $\beta (C_\varphi) = \beta(C_\psi)$, 
and therefore, with the notation of \eqref{Dieze}: 
\begin{displaymath}
\qquad \quad \beta (C_\varphi) \geq [\varphi^{\#} (a)]^2 , \quad \text{for all } a \in \D.
\end{displaymath}
By passing to the supremum on $a\in\D$,  we end the proof of   Proposition~\ref{Canon}, and that of Theorem~\ref{Secondary} in the compact case. If 
$C_\varphi$ is not compact, the proposition trivially holds. Indeed, in this case, we have $\beta (C_\varphi) = 1 \geq [\varphi]^2$. \qed
\bigskip

\noindent{\bf Remark.} It is easy to see that the composition operator $C_\varphi$ is always of infinite rank, contrary to the case of a Hankel operator, 
so that in some sense it refuses to be approached by finite-rank operators. Theorem~\ref{Secondary} quantifies things: it is a well-known and easy fact  
(see for example \cite{Shap-livre}, page~25 and see Theorem~\ref{Ternary} to come) that, in the case $\Vert \varphi \Vert_\infty < 1$, we have 
$a_{n} (C_\varphi) \leq c \Vert \varphi \Vert_\infty^n$ (and hence $\beta (C_\varphi) \leq \| \varphi \|_\infty < 1$), showing that the approximation 
numbers can decrease at an exponential speed. Theorem~\ref{Secondary} shows that this speed is the maximal possible one. The next theorem says that this 
maximal speed is \emph{only} obtained when  $\Vert \varphi \Vert_\infty < 1$. 
\par

\begin{theoreme}\label{Lastresult} 
For every $\alpha \geq - 1$, there exists, for any $0 < r < 1$, $s = s(r) < 1$, satisfying  $\lim_{r\to 1^-} s (r) = 1$, such that, for   
$C_\varphi \colon {\mathfrak{B}_\alpha}\to {\mathfrak{B}_\alpha}$, one has, with the notation coined in \eqref{beta}: 
\begin{equation} 
\Vert \varphi \Vert_\infty > r \quad \Longrightarrow \quad \beta (C_\varphi)\geq s^2 .
\end{equation} 
In particular, the exponential speed of convergence to $0$ of the approximation numbers of a composition operator $C_\phi$ of symbol $\varphi$ takes place 
if and only if $\Vert \varphi\Vert_\infty < 1$; in other words, we have:
\begin{equation}\label{Onlycase} 
\Vert \varphi\Vert_\infty = 1 \quad \Longleftrightarrow  \quad \beta (C_\varphi) = 1.
\end{equation}
\end{theoreme} 
\smallskip

The proof will proceed through a series of lemmas. Throughout that proof, we assume, without loss of generality, that $\varphi (0) = 0$. 
\goodbreak

\begin{lemme} \label{Seville 1} 
Let $K$ be a compact subset of $\varphi (\D)$ and $\mu$ be a probability supported by $K$. Then, there exists a constant $\delta > 0$ such that, if  
$R_\mu \colon \mathfrak{B}_\alpha \to L^{2}(\mu)$ denotes the restriction operator, we have:
\begin{displaymath} 
a_{n} (C_\varphi) \geq \delta \, a_{n} (R_\mu). 
\end{displaymath} 
In particular:
\begin{displaymath} 
\beta (C_\varphi) \geq \beta (R_\mu).
\end{displaymath} 
\end{lemme}

\noindent{\bf Proof.} Since $\varphi$ is an open map, there exists a compact set $L\subset \D$ and a Borel subset $A \subset L$ such that $\varphi (A) = K$ 
and $\varphi \colon A \to K$ is a bijection (see \cite{Pa}, Chapter I, Theorem~4.2). Then $\mu = \varphi (\nu)$, where $\nu = \varphi^{-1} (\mu)$ is a 
probability measure supported by $L$, and we have automatically $\Vert R_\nu \Vert < \infty$. Then, for every $f\in {\mathfrak B}_\alpha$:
\begin{displaymath} 
\Vert f \Vert_{L^{2}(\mu)}^2 = \int_{K} \vert f \vert^2 \, d\mu = \int_{L} \vert f \circ \varphi \vert^2 \, d\nu 
= \Vert C_{\varphi}f \Vert _{L^{2}(\nu)}^2.
\end{displaymath} 
This yields $\Vert R_{\mu}f \Vert = \Vert (R_{\nu} \circ C_{\varphi}) f \Vert$, so $C_\varphi$ acts as an isometry from $L^2 (\mu)$ into $L^2 (\nu)$, and 
the lemma follows, since we have then: 
\begin{displaymath} 
a_{n}(R_\mu) = a_{n}(R_\nu \circ C_\varphi) \leq \Vert R_\nu \Vert \, a_{n} (C_\varphi) 
\end{displaymath} 
for every $n\geq 1$. \hfill $\square$
\par
\medskip

Observe that this provides a new proof of Theorem~\ref{Secondary}. Indeed, if $K \subset \varphi (\D)$ is a small ball of center $0$ and radius $r$, we can take 
for $\mu$ the normalized area measure on $K$; then Parseval's formula easily shows that $\beta (R_\mu) \geq r$ in that case. \par
\smallskip

The strategy of the proof of Theorem \ref{Lastresult} will consist of refining this observation. More precisely, we shall show that  the situation can be reduced to 
the case $K = [0, r]$, and that  an appropriate  choice of  $\mu$ can be made in that case, giving a sharp lower bound for $\beta (R_\mu)$. We begin with 
explaining that choice in the next two lemmas.

\begin{lemme} \label{Seville 2} 
For every $r  \in (0, 1)$ there exists $s = s (r) < 1$ and $f = f_r \in H^\infty$ with the following properties:\par \smallskip 
1) $\lim_{r \to 1^-}s(r) = 1$; \par \smallskip 
2) $\Vert f \Vert_\infty \leq 1$; \par \smallskip 
3) $f ((0, r]) = s\, \partial \D$ in a one-to-one way.
\end{lemme}

\noindent{\bf Proof.} Let $\rho = \frac{1 - \sqrt{1 - r^2}}{r}$\,. Then $r = \frac{2 \rho}{1 + \rho^2}$ and the automorphism 
$\varphi_{\rho} (z) = \frac{\rho - z}{1 - \rho z}$ maps $[0, r]$ onto $[- \rho, \rho]$. We define $\varepsilon =\varepsilon (r)$ and $s = s(r)$ by the 
following relations: 
\begin{equation} 
\varepsilon (r) = \frac{\pi}{\log \frac{1 + \rho}{1 - \rho}} \,\raise 1pt \hbox{,} \qquad \text{and} \quad s = \e^{-\varepsilon \pi / 2}.
\end{equation} 
Let now 
\begin{equation} 
\chi (z) = \varepsilon \log \frac{1 + \varphi_{\rho}(z)}{1 - \varphi_{\rho}(z)} 
\end{equation} 
and 
\begin{equation} 
f (z) =  s\, \e^{i \chi(z)}\,.
\end{equation} 
Note that $f = \e^h$, where 
\begin{displaymath} 
h (z) = i\varepsilon \log  \frac{1 + \varphi_{\rho} (z)}{1 - \varphi_{\rho} (z)} - \varepsilon \, \frac{\pi}{2}
\end{displaymath} 
is a conformal mapping from $\D$ onto a small vertical strip of the left-half plane. This function $f$ fulfills all the requirements of the lemma. Indeed, we have 
$\vert f (z) \vert \leq 1$ for all $z\in\D$ and
\begin{displaymath} 
h ([0, r]) = \bigg [- i \varepsilon \log \frac{1 + \rho}{1 - \rho}\, \raise 1pt \hbox{,} 
\ i \varepsilon \log \frac{1 + \rho}{1 - \rho}\bigg]  - \varepsilon\, \frac{\pi}{2} 
= [- i \pi , i \pi] - \varepsilon\, \frac{\pi}{2}\, \raise 1pt \hbox{,} 
\end{displaymath} 
so that $f ((0, r]) = \{w = s \e^{i\theta}\,;\  -\pi \leq \theta \leq \pi\}$, in a one-to-one way. \hfill $\square$
\par
\bigskip

Lemma~\ref{Seville 2} allows a good choice of the measure $\mu$ as follows.

\begin{lemme}\label{Seville 3} 
Let $f$ be as in  Lemma~\ref{Seville 2}. Then, there  exists a probability measure $\mu =\mu_r$ supported by $[0, r]$ and a constant $\delta_r > 0$ such that, 
for any integer $n \geq 1$ and any choice of scalars $c_0, c_1, \ldots, c_{n-1}$, we have:
\begin{displaymath} 
\bigg\Vert \sum_{j = 0}^{n -1}c_j R_{\mu}(f^j) \bigg\Vert _{L^{2}(\mu)} 
\geq  \frac{s^n}{\sqrt n} \bigg\Vert \sum_{j = 0}^{n - 1}c_j f^j \bigg\Vert _{H^2} 
\geq  \frac{s^n}{\sqrt n} \bigg\Vert \sum_{j = 0}^{n - 1}c_j f^j \bigg\Vert _{{\mathfrak B}_\alpha} .
\end{displaymath} 
As a consequence, we can claim that, for $C_\varphi \colon {\mathfrak B}_\alpha \to {\mathfrak B}_\alpha$:
\begin{equation} \label {Firstfact} 
\varphi (\D) \supset [0, r] \quad \Longrightarrow  \quad \beta (C_\varphi) \geq s = s (r) .
\end{equation}
\end{lemme}

\noindent{\bf Proof.} With our previous notations, we know that $\chi$ is a bijective map from $]0, r]$ onto the unit circle $\partial \D$. Let 
$\mu = \chi^{-1} (m)$ be the image of the Haar measure $m$ of $\partial \D$ by $\chi^{-1}$. We have by definition of $\mu$: 
\begin{align*}
\bigg\Vert \sum_{j = 0}^{n - 1}c_j R_{\mu}(f^j) \bigg\Vert _{L^{2}(\mu)}^2 
& =\int_{0}^r \bigg\vert \sum_{j = 0}^{n - 1}c_j f^{j} (x) \bigg\vert^2 \, d\mu (x) 
=\int_{0}^r \bigg\vert \sum_{j = 0}^{n - 1}c_j s^j  \e^{i j \chi (x))} \bigg\vert^2 \,d\mu (x) \\ 
& = \int_{-\pi}^\pi \bigg\vert \sum_{j = 0}^{n - 1}c_j s^j \e^{i j \theta} \bigg\vert ^{2} \, \frac{d\theta }{2\pi} 
=\sum_{j = 0}^{n - 1} \vert c_j \vert^2 s^{2j} \\
& \geq  s^{2n} \sum_{j = 0}^{n - 1}\vert c_j \vert^2. 
\end{align*}
Now, $\| f^j\|_{H^2} \leq \| f ^j \|_\infty \leq 1$, so that we have, using the Cauchy-Schwarz inequality:
\begin{displaymath} 
\bigg\Vert \sum_{j = 0}^{n - 1}c_j f^j \bigg\Vert _{H^2} \leq \sum_{j = 0}^{n - 1} \vert c_j \vert \, \Vert f^j \Vert_{H^2} 
\leq \sum_{j = 0}^{n - 1} \vert c_j \vert 
\leq \sqrt n \bigg(\sum_{j = 0}^{n - 1} \vert c_j \vert^2 \bigg)^{1/2} ,  
\end{displaymath} 
giving the first inequality, since $\Vert \ \Vert_{H^2} \geq \Vert \ \Vert_{\mathfrak{B}_\alpha}$. Finally, let 
$R \colon {\mathfrak B}_\alpha \to L^{2}(\mu)$ be an operator of rank $< n$. We can find a function 
$g = \sum_{j = 0}^{n - 1}c_j f^j$ such that $\Vert g \Vert_{{\mathfrak B}_\alpha} = 1$ and $R (g) = 0$. The first part of the proof gives:
\begin{align*}
\Vert R_\mu - R \Vert 
& \geq \Vert R_\mu (g) - R (g) \Vert = \Vert R_\mu (g) \Vert = \bigg\Vert \sum_{j = 0}^{n - 1} c_j f^j \bigg\Vert_{L^{2}(\mu)} \\
& \geq \frac{s^n}{\sqrt n}\bigg\Vert \sum_{j = 0}^{n - 1} c_{j} f^j \bigg\Vert_{{\mathfrak B}_\alpha} 
= \frac{s^n}{\sqrt n} \,\cdot
\end{align*}
Therefore $a_{n}(R_\mu) \geq s^n /\sqrt n$ and, in view of Lemma \ref{Seville 1}, the last conclusion of Lemma \ref{Seville 3} 
follows. \hfill $\square$  
\par\bigskip

The next lemma explains how to reduce the situation to the case $K = [0, r]$ when we only know that $\Vert \varphi \Vert_\infty > r$. It was inspired to us by 
the proof of the Lindel\"of theorem that convergence along a curve implies  non-tangential  convergence for functions in Hardy spaces 
(\cite{Rudin-livre} page 300). 

\begin{lemme} \label{Seville 4} 
Suppose that $0$ and $r$ belong to $\varphi (\D)$, with $0 < r < 1$.  Let $\mu$ be a probability measure carried by $[0,r]$.  Then, there exists a probability measure 
$\nu$ carried by a compact set $K \subset \varphi (\D)$ such that, for any $f \in \mathcal{H} (\D)$:
\begin{equation} \label{Harmonic} 
\int_{[0, r]} \vert f (x) \vert^{2} \,d\mu(x) \leq \frac{1}{2} \int_{K} \big( \vert f (z) \vert^2 + \vert f (\bar{z}) \vert^2 \big) \,d\nu (z) . 
\end{equation}
\end{lemme}

\noindent{\bf Proof.} Since $\varphi (\D)$ is open and connected and $0, r \in \varphi (\D)$, there is a curve with image $K \subset \varphi(\D)$ connecting 
$0$ and $r$. Put $\tilde K = \{\bar{z} \, ; \  z \in K\}$. Then, there exists a compact set $L$ such that $[0, r] \subset L$ and whose boundary 
$\partial L \subset (K \cup \tilde K)$. Now, the existence of $\nu$ carried by $K$ will be provided by an appropriate application of the Pietsch factorization 
Theorem. To that effect,  let $X$ be the real subspace of $\mathcal{C} (L)$ formed by the real functions which are harmonic in the interior of $L$. By the 
maximum principle for harmonic functions, $X$ can be viewed as a subspace of $\mathcal{C} (K \cup \tilde K)$. Now,  the inclusion map $j$ of $X$ into 
$L^{2}(\mu)$ has $2$-summing norm less than one (\cite{AL-KA-livre} page 208, or \cite{Li-Queff}, Chapitre 5, Proposition~I.3). Therefore, the Pietsch 
factorization Theorem (\cite{AL-KA-livre} page 209, or \cite{Li-Queff}, Chapitre 5, Th\'eor\`eme~I.5) implies the existence of a probability $\sigma$ on 
$K\cup\tilde K$ such that, for every $u \in X$: 
\begin{equation}\label{One} 
\Vert u \Vert_{L^{2} (\mu)}^2 = \int_{[0, r]} u^{2} \, d\mu 
\leq \int_{K \cup \tilde K} u^{2} \,d\sigma .  
\end{equation}
For any harmonic function $u$ on $\D$, we can apply \eqref{One} to $u (z)$ and $u (\bar{z})$ to get:
\begin{displaymath} 
2 \int_{[0, r]} u^{2} \, d\mu 
\leq  \int_{K \cup \tilde K} \big[ u^{2} (z) + u^{2} (\bar{z}) \big] \, d\sigma (z) 
= \int_{K \cup \tilde K} \big[ u^{2} (z) + u^{2} (\bar{z}) \big] \, d\tilde\sigma (z) ,
\end{displaymath} 
where $\tilde \sigma$ is the symmetric measure of $\sigma$, defined by $\tilde \sigma(E) = \sigma(\bar{E})$. There is a probability $\nu$ on $K$ such that 
$\nu + \tilde \nu = \sigma + \tilde \sigma$. For this probability $\nu$, we thus have, for any real  harmonic function $u$ on $\D$:
\begin{equation}\label{Two} 
\Vert u \Vert_{L^{2} (\mu)}^2 
\leq  \int_{K} \big[ u^{2} (z) + u^{2} (\bar{z}) \big] \, d\nu (z).   
\end{equation}
Now, given $f \in \mathcal{H} (\D)$, we use \eqref{Two} with $u$ the real and imaginary parts of $f$, and sum up to get \eqref{Harmonic}. 
\hfill $\square$
\bigskip

We can now finish the proof of Theorem~\ref{Lastresult} as follows. \par 
Suppose that $\Vert \varphi \Vert_\infty > r$. Then, making a rotation if necessary, 
we may assume that $0, r \in \varphi (\D)$ (recall that $\varphi (0) = 0$). Let $\mu$ as in Lemma~\ref{Seville 3}. Using Lemma~\ref{Seville 4}, we find a 
probability measure $\nu$, compactly supported by $\varphi (\D)$, such that \eqref{Harmonic} holds. This inequality shows that:
\begin{displaymath} 
\Vert R_{\mu} f \Vert ^2 \leq \frac{1}{2} \big( \Vert R_{\nu} f \Vert ^2 + \Vert R_{\tilde\nu} f \Vert ^2 \big),
\end{displaymath} 
so that $R_\mu = A (R_\nu \oplus R_{\tilde\nu})$ with $\Vert A \Vert \leq 1/\sqrt 2 \leq 1$. Therefore, by the ideal and sub-additivity properties 
\eqref{ideal}:
\begin{displaymath} 
a_{2n} (R_\mu) \leq a_{2n} (R_\nu \oplus R_{\tilde\nu}) \leq a_{n} (R_\nu) + a_{n} (R_{\tilde\nu}) = 2 \, a_{n} (R_\nu) ,
\end{displaymath} 
implying $\beta (R_{\nu}) \geq \beta (R_{\mu})^2$.  Finally, Lemma~\ref{Seville 1} and Lemma~\ref{Seville 3} give:
\begin{displaymath} 
\beta (C_\varphi) \geq \beta (R_{\nu}) \geq \beta (R_{\mu})^2  \geq s (r)^2 , 
\end{displaymath} 
and this ends the proof of Theorem~\ref{Lastresult}.  \hfill $\square$
\bigskip

\noindent{\bf Remark.} The proof of Theorem~\ref{Lastresult} is strongly influenced by the papers \cite{ERO} and \cite{WID}. In the first one, it is proved 
that, if $K$ is a continuum of a connected open set $\Omega$ and if the doubly connected region $\Omega \setminus K$ is conformally equivalent to the 
annulus $1 < \vert z \vert < R$, then  there exists a linearly independent sequence $(f_n)$ in $H^{\infty}(\Omega)$ satisfying, for all scalars $c_j$: 
\begin{displaymath} 
\bigg\Vert \sum_{j = 1}^n c_j f_j \bigg\Vert_{H^{\infty}(\Omega)} \leq R^n \, \bigg\Vert \sum_{j = 1}^n c_j f_j \bigg\Vert_{\mathcal{C}(K)}.
\end{displaymath} 
As a consequence, the author proves that $\lim_{n\to \infty} d_{n}^{1/n}= 1/ R$, where the numbers $d_n$ are the Kolmogorov numbers 
(see \cite{CA-ST-livre} page~49 for the definition) of the restriction map $H^{\infty}(\Omega) \to \mathcal{C}(K)$. This statement led us to 
Lemma~\ref{Seville 3}. In the second paper, it is proved that, for the same operator, one has 
$\lim_{n\to \infty} d_{n}^{1/n} = \e^{-1 / C (K, \Omega)}$, where $C (K, \Omega)$ is the Green capacity of $K$ relative to $\Omega$. So that one has 
$1/ R = \e^{-1 / C (K, \Omega)}$. In the case we were interested in, namely $\Omega = \D$ and $K_r = [- r, r]$, it seemed to us, for topological and 
analytic reasons, that $R$ should tend to $1$ as $r \to 1$, in other terms that we should have  $\lim_{r \to 1^-} C (K_r, \D) = \infty$. This is indeed the case 
(\cite{Saff-Totik-livre}, Example~II.1), but the proof is fairly involved, and the desire to get a reasonably simple and self-contained proof of 
Theorem~\ref{Lastresult} led us to the previous series of lemmas, once we were sure that the result was true. 

\goodbreak

\section{Slow speed} \label{section 4} 

In this section, we shall see that the speed of convergence to $0$ of the approximation numbers of a compact composition operator can be as slow as one wants. 
This answers in the positive to a  question which was first asked to us by C.~Le~Merdy (\cite{Le Merdy}) in the OT Conference 2008 of Timisoara. 

\begin{theoreme} \label{Principal} 
Let $(\varepsilon_n)_{n\geq 1}$ be a non-increasing sequence of positive real numbers of limit zero. Then, there exists an injective Schur function $\varphi$ 
such that $\varphi(0) = 0$ and  $C_{\varphi} \colon \mathfrak{B}_\alpha\to\mathfrak{B}_\alpha$ is compact, i.e. $a_{n} (C_{\varphi})\to 0$,  but:
\begin{equation} \label{Main} 
\liminf_{n \to \infty} \frac{a_{n} (C_\varphi)}{\varepsilon_{n}} > 0. 
\end{equation}
Equivalently, we  have for some  positive number $\delta>0$, independent of $n$:
\begin{displaymath}
\qquad \qquad \qquad a_{n} (C_\varphi) \geq \delta \, \varepsilon_{n} \qquad  \text{for all } n\geq 1.
\end{displaymath}
\end{theoreme}
\medskip

As in the case of Hankel operators, an immediate consequence of  Theorem~\ref{Principal} is the following:
\begin{corollaire} \label{Coco} 
There exists a  composition operator $C_\varphi \colon H^2\to H^2$ which is compact, but in no Schatten class.
\end{corollaire}
This corollary, which Theorem~\ref{Principal} reinforces and precises, was an answer to a question of Sarason, and has been first proved in  
\cite{Carroll-Cowen}. Other proofs appeared in \cite{A}, \cite{JO}, \cite{JFA}, \cite{JMAA}, \cite{Zhuu} (for a positive result on Schatten-ness, we refer to   
\cite{LuZh}).
\medskip

The construction of  the symbol $\varphi$ in Theorem~\ref{Principal} follows that given in \cite{Carroll-Cowen}, but we have to proceed to some necessary 
adjustments. In order to exploit \eqref{Hyperbolic}, we shall use, as in \cite{Carroll-Cowen}, the following two results due to Hayman (\cite{Hayman}) 
concerning the hyperbolic distance $d (z, w; U)$ of two points $z, w$ of a simply connected domain $U$ (see also \cite{Keen-livre}), whose proof uses in 
particular the comparison principle~\eqref{Principle}:  
\begin{proposition}\label{Majorant} 
Suppose that $U$ contains the rectangle 
\begin{displaymath}
R = \{z \in \C\, ;\ a_1 - b < \Re z < a_2 + b,\, \vert \Im z \vert < b\} ,
\end{displaymath}
where $a_1 < a_2$ and $b > 0$. Then, we have the upper estimate:
\begin{equation} \label{Upper} 
d (a_1, a_2; U)\leq \frac{\pi}{4b} (a_2 - a_1) + \frac{\pi}{2} \cdot 
\end{equation}
\end{proposition}

\begin{proposition}\label{Minorant} 
Suppose that   $U$ contains the rectangle 
\begin{displaymath}
R = \{z \in \C\, ;\ a_1 - c < \Re z < a_2 + c, \, \vert \Im z \vert < c\}, 
\end{displaymath}
where $a_1 < a_2$ and $c > 0$, but that the horizontal sides 
\begin{displaymath}
\{z \in \C \, ; \ a_1 - c \leq \Re z \leq a_2 + c, \, \vert \Im z \vert = c\} 
\end{displaymath}
of that rectangle are disjoint from $U$. Then, we have the lower estimate:
\begin{equation} \label{Lower} 
d (a_1, a_2 ; U) \geq \frac{\pi}{4c} (a_2 - a_1) - \frac{\pi}{2} \cdot
\end{equation}
\end{proposition}
\medskip

We now proceed to the construction of our Schur function $\varphi$. 
\par

We first define a continuous  map $\psi \colon\R\to \R$ as follows: 
\begin{displaymath}
\psi (t) = 
\begin{cases}
K (1 + \vert t \vert) & \mathrm{if}\ \vert t \vert \leq 1 \\ 
\vert t \vert / A (\vert t \vert) &\mathrm{if}\ \vert t \vert > 1,
\end{cases}
\end{displaymath}
where $K$ is a positive constant adjusted below and $A \colon [0,\infty[ \to [0,\infty[$ an increasing piecewise linear function on the intervals $(0,1)$ 
and $(\e^{n - 1}, \e^n)$ such that  
\begin{displaymath}
A (0) \mathop{=}^{def} A_0 = 0, \quad A (\e^{n - 1}) \mathop{=}^{def} A_n  \quad \text{for }  n \geq 1,  \quad \text{and}  \quad 2K = 1/ A(1), 
\end{displaymath}
the increasing sequence $(A_n)$ being positive and concave for $n\geq 1$, and tending to $\infty$. It then follows that the sequence of slopes 
$\frac{A_n - A_{n - 1}} {\e^n - \e^{n - 1}}$ is decreasing, since $A_{n+1} - A_{n} \leq A_{n} - A_{n-1} \leq \e \, (A_{n} - A_{n - 1})$, that the function 
$A$ is increasing and concave on $(0,\infty)$ and vanishing at $0$, implying that $A(t)/ t$ is decreasing on $(0,\infty)$, and that in particular  $\psi$ is 
increasing on $(1, \infty)$.\par

We then define a domain $\Omega$ of the complex plane by:
\begin{equation} \label{Domain} 
\Omega = \{w\in\C \, ; \ \vert \Im w \vert < \psi (\vert \Re w \vert)\}.
\end{equation}
Let $\sigma \colon \D \to \Omega$ be the unique Riemann map such that $\sigma(0) = 0$ and $\sigma' (0) > 0$. This map exists in view of the following 
simple fact.

\begin{lemme} \label{Star} 
The domain  $\Omega$ defined by \eqref{Domain} is star-shaped with respect to the origin and $\sigma \colon (-1, 1) \to \R$ is an increasing bijection 
such that $\sigma (-1) = -\infty$ and $\sigma (1) = \infty$.
\end{lemme}

\noindent{\bf Proof.} The star-shaped character of $\Omega$ will follow from the implication: 
\begin{displaymath}
\vert \Im w \vert < \psi (\vert \Re w \vert) \quad \text{and} \quad 0 < \lambda < 1 
\quad \Longrightarrow \quad \vert \Im (\lambda w) \vert < \psi (\vert \Re (\lambda w) \vert).
\end{displaymath}
We may assume that both $\Re w, \Im w$ are positive, and it is enough to prove:
\begin{equation} \label{Shaped} 
\lambda \psi(x) \leq \psi (\lambda x) ,  \qquad \forall \lambda \in [0,1], \ \forall x > 0.
\end{equation}
This is easy to check separating three cases:\par\smallskip
1)  $x\leq 1$; then $\lambda \psi (x) = \lambda K (1 + x) \leq K (1 + \lambda x) = \psi (\lambda x)$;\par\smallskip
2)  $\lambda x \leq 1 < x$; then, since $A (x) > A (1)$, 
\begin{displaymath}
\lambda \psi (x) = \lambda \frac{x}{A(x)} < 2K \lambda x \leq K (1+ \lambda x) = \psi (\lambda x);
\end{displaymath}
\par 
3)  $\lambda x > 1$; we then have, since $\psi$ increases, 
\begin{displaymath}
\lambda \psi (x) = \lambda \frac{x}{A(x)} \leq \frac{\lambda x}{A (\lambda x)} = \psi (\lambda x)
\end{displaymath}
and  this ends the proof of \eqref{Shaped}. Now, since $\sigma$ is determined by the value of  $\sigma (0)$ and  the sign of $\sigma' (0)$, we have 
$\sigma (\overline{z}) = \overline{\sigma (z)}$ for all $z\in\D$, so that $\sigma [(-1, 1)] \subset \R$. And since the derivative of an injective analytic 
function does not vanish and $\sigma' (0) > 0$, we get that $\sigma$ is increasing on $(-1, 1)$. Finally, if $w \in \R$ and $w = \sigma (z)$, we have 
$\overline{w} = w$, so that $\sigma (\overline{z}) =\sigma (z)$ and $\overline{z} = z$, which proves the surjectivity of $\sigma \colon (-1,1) \to \R$.
\qed 
\medskip

We now choose  $A_n$ as follows, $\eta > 0$ denoting a positive numerical constant to be specified later.
\begin{equation} \label{Choice} 
\qquad A_n = \eta \log \frac{1}{\varepsilon_n} \,\raise 1,5pt \hbox{,} \quad  n\geq 1.
\end{equation}
Observe that this is an increasing, concave sequence tending to $\infty$ since we assumed that $(\varepsilon_n)$ is log-convex and decreasing to $0$.
\par
Finally,  we define our Schur function $\varphi$ and our sequence $(r_n)$ under the form of the following lemma, in which the increasing character of $\psi$ 
is important.
\begin{lemme} \label{Summarize} 
Let $\varphi$ be defined by  
\begin{displaymath}
\varphi (z) = \sigma^{-1} (\e^{-1} \sigma (z)),  
\end{displaymath}
and let $r_n = \sigma^{-1} (\e^n)$. Then we have: \par \smallskip 
{\rm 1.} $\varphi$ is univalent and maps $\D$ to $\D$, $(r_n)$ increases, and $\varphi (0) = 0$; \par\smallskip
{\rm 2.} $\varphi (r_{n+1}) = r_n$;\par\smallskip 
{\rm 3.} $\frac{1 - r_{n+1}}{1 - r_n} \to 0$ and therefore $(r_n)$ is an interpolation sequence;\par\smallskip
{\rm 4.} $C_\varphi \colon \mathfrak{B}_\alpha \to \mathfrak{B}_\alpha$ is compact. 
\end{lemme}
\noindent{\bf Proof.} 1. Since $\Omega$ is star-shaped, $\e^{-1}\sigma (z) \in \Omega$ when $z \in \D$, so $\varphi$ is well-defined and maps 
$\D$ to itself in a univalent way. Moreover, $\varphi (0) = \sigma^{-1} (0) = 0$, and $(r_n)$ increases since $\sigma^{-1}$ increases on $\R$.\par
\smallskip
2. We have 
$\varphi (r_{n+1}) = \sigma^{-1} \big(\frac{1}{\e} \sigma (r_{n+1}) \big) = \sigma^{-1} \big(\frac{1}{\e} \,\e^{n+1} \big) 
= \sigma^{-1} (\e^{n}) = r_n$.
\par\smallskip
3. This assertion  is more delicate and relies on Proposition~\ref{Minorant} as follows.\par 

Set $d_n = \psi (\e^n)$. We have clearly  $\e^{n+1} + d_{n+2} < \e^{n+2}$ for large $n$ (recall that $\psi (t) = o\,(t)$ as $t \to \infty$), so that 
$\psi (\e^{n+1} + d_{n+2}) < \psi (\e^{n+2}) = d_{n+2}$ since  $\psi$ is increasing. By the intermediate value theorem for the function 
$\psi (\e^{n+1} + x) - x$, we can therefore find a positive number   $c_n < d_{n+2}$  such that $\psi (\e^{n+1} + c_n) = c_n$. \par

Now, consider  the open sets:
\begin{displaymath}
R_{n} = \{z \in \C \, ;\ \e^n - c_n < \Re z < \e^{n+1} + c_n  \text{ and } \vert \Im z \vert < c_n\}, \quad U_{n} = R_{n} \cup \Omega .
\end{displaymath} 
Those sets $U_{n}$ satisfy the assumptions of Proposition~\ref{Minorant} in view of \eqref{Domain}. Indeed, if $z$ belongs to the horizontal sides of $R_n$, 
we have $z \notin U_n$ since
\begin{displaymath}
\e^n - c_n \leq \Re z \leq \e^{n+1} + c_n \quad \Longrightarrow \quad \psi (\Re z) \leq \psi (\e^{n+1} + c_n) = c_n =\vert \Im z \vert .
\end{displaymath}
This proposition then  gives, since $\Omega \subset U_{n}$ and $c_n < d_{n+2}$, and since the hyperbolic metric is conformally invariant,
\begin{align*}
 d (r_n, r_{n+1}) 
& = d (\e^n, \e^{n+1}; \Omega) \geq d (\e^n, \e^{n+1}; U_{n}) \geq \frac{\pi}{4c_n} (\e^{n+1} - \e^n) - \frac{\pi}{2} \\
& \geq c\, \frac{\e^{n+2}}{\psi (\e^{n+2})} = c A (\e^{n+2}) \geq c A_n, 
\end{align*}
where $c$ is a positive constant. Now, we use Lemma~\ref{Poincare} to obtain:
\begin{displaymath}
\frac{1 - r_{n+1}}{1 - r_{n}} \leq 2 \, \e^{-2 d (r_n,  r_{n+1})} \leq 2 \,\e^{-2c A_n},
\end{displaymath}  
which proves that $\frac{1 - r_{n+1}}{1 - r_{n}} \to 0$, and implies that $(r_n)$ is an interpolation sequence.
\par\smallskip

4. Since $\varphi$ is univalent, the compactness of $C_\varphi \colon \mathfrak{B}_\alpha \to \mathfrak{B}_\alpha$  amounts to proving that 
$\lim_{\vert z \vert \to 1} \frac{1 - \vert \varphi (z) \vert} {1 - \vert z \vert} =\infty$. For $\alpha > -1$, this follows from \cite{MASH}, Theorem~3.5 
and for $\alpha = -1$ from \cite{Shap-livre}, page~39.  By the Julia-Carath\'eodory Theorem (\cite{Shap-livre}, page~57), this in turn is equivalent to proving 
that  for any $u, v$ on the unit circle, the quotient $\frac{\varphi (z) - v}{z - u}$ has no finite limit as $z$ tends to $u$ radially. This latter fact requires some 
precise justification.\par

First, we notice that $\sigma$ extends continuously to an injective map of the open upper half of the unit circle onto the upper part of the boundary of 
$\Omega$ (and similarly for lower parts). This follows from the Carath\'eodory extension theorem (\cite{Rudin-livre}, page~290), applied to the restriction of 
$\sigma^{-1}$ to the Jordan region limited by $\partial \Omega$ and two vertical lines $\Re w = \pm R$ where $R > 0$ is arbitrarily large. Now, let 
$u \in \partial \D$ with $u \neq \pm 1$. Then, $\sigma (ru) \to w \in \partial \Omega$ as $r \to 1^-$, so that 
$\e^{-1} \sigma (ru) \to \e^{-1} w = w' \in \Omega$ and that $\varphi (ru) \to \sigma^{-1} (w') \in\D$. Therefore the image of $\varphi$ touches the 
unit circle only at $\pm 1$, and the assumption of the Julia-Carath\'eodory Theorem is fulfilled if  $u \neq \pm 1$. By symmetry, it remains to test the point 
$u = 1$ for which we have:
\begin{displaymath}
\limsup_{r \mathop{\longrightarrow}\limits^{<} 1}   \frac{1 -\varphi (r)}{1 - r} 
\geq \limsup_{n\to\infty} \frac{1 - \varphi (r_{n+1})}{1 - r_{n+1}} 
= \limsup_{n\to\infty} \frac{1 - r_{n}}{1 - r_{n+1}} = \infty 
\end{displaymath}
by the preceding  point~3. Since $\vert v - \varphi (r) \vert \geq 1 - \varphi (r)$, this ends the proof of Lemma~\ref{Summarize}.\qed
\bigskip

We now want a good lower bound for the weights $w_n$ appearing in \eqref{Good Weights}. To that effect, we apply Proposition~\ref{Majorant} with 
\begin{displaymath}
U = \Omega, \quad a_1 = \e^n, \ a_2 = \e^{n+1} \quad \text{and} \quad b_n = \psi (\e^{n - 1}), 
\end{displaymath}
as well as
\begin{displaymath}
R'_{n} = \{z \in \C \,; \ \e^n - b_n < \Re z < \e^{n+1} + b_n  \text{ and } \vert \Im z \vert < b_n\} .
\end{displaymath}
We have $\e^n - b_n > \e^{n - 1}$ for large $n$, since this amounts to  
\begin{displaymath}
\e^n - \e^{n - 1} >  b_n = \frac{\e^{n - 1}}{A (\e^{n - 1})} \, \raise 1,5pt \hbox{,}  \quad 
\text{or} \quad \e - 1 > \frac{1}{A (\e^{n - 1})} \raise 1,5pt \hbox{,} 
\end{displaymath}
which holds for large $n$ since $A (t)$ tends to $\infty$ with $t$. We then observe that $R'_{n} \subset \Omega$. Indeed, 
$z \in R'_{n} \Longrightarrow \Re z >  \e^n  - b_n > \e^{n - 1}$ and, since $\psi$ is increasing, we have 
$\psi (\Re z) > \psi (\e^{n - 1}) = b_n > \vert \Im z \vert$. Therefore, we can apply \eqref{Upper} and get, for all $n \geq 1$: 
\begin{displaymath}
d (\e^n, \e^{n + 1} ; \Omega) \leq \frac{\pi}{4 \psi (\e^{n - 1})} (\e^{n+1} - e^n) + \frac{\pi}{2} \leq C_0 A (\e^{n - 1}) = C_0 A_n  , 
\end{displaymath}
where $C_0$ is a \emph{numerical} constant. By conformal invariance, we have as well $d (r_n,  r_{n+1}) \leq C_0 A_n$. It then follows from 
\eqref{Hyperbolic} that: 
\begin{equation} \label{Hayman} 
\frac{1 - r_{n + 1}} {1 - r_{n}} \geq \exp \big(-2 d (r_n, r_{n+1}) \big) \geq \exp (-2 C_0 A_n). 
\end{equation}

Now, we take $h(z) = z - r_1$ in Lemma~\ref{Simple minoration} and use the  ideal property \eqref{ideal} of the approximation numbers. We get, denoting 
by  $C$ the interpolation constant of the sequence $(r_n)$, and using the fact that $\Vert M_h \Vert = \Vert h \Vert_\infty \leq 2$:
\begin{equation} \label{Comparison} 
a_{n} ({\mathbf B}) \leq \Vert J^{-1} \Vert \, a_{n} (C_{\varphi}) \,\Vert M_h\Vert \, \Vert J \Vert \leq 2 C^2 a_{n} (C_{\varphi}).
\end{equation} 
Next, we choose $\eta = 1/ C_0$ in \eqref{Choice} and we set  $d = (r_2 - r_1)/\sqrt 2$. Using  Lemma~\ref{Compression} and relations 
\eqref{Bergman}, \eqref{Good Weights} and \eqref{Hayman}, we see  that  the weights $w_n$ associated with ${\mathbf B}$ verify:  
\begin{equation} \label{Weights} 
\begin{split}
\vert w_n \vert 
& = h (r_{n + 1}) \frac{\Vert K_{r_n}\Vert}{\Vert K_{r_{n + 1}} \Vert} 
= h (r_{n+1}) \sqrt{\frac{1 - r_{n + 1}^2} {1 - r_{n}^2}} \geq \frac{r_2 - r_1}{\sqrt 2} \sqrt{\frac{1 - r_{n + 1}} {1 - r_n}} \\
& \geq d \exp (- C_0 A_n) \geq d \varepsilon_n \qquad \text{for\ all }  n \geq 1.
\end{split} 
\end{equation}
Finally, using Lemma~\ref{Simple minoration}, \eqref{Comparison} and \eqref{Weights}:
\begin{displaymath}
a_{n} (C_{\varphi}) \geq \frac{1}{2 C^2} a_{n} ({\mathbf B}) \geq \frac{1}{2 C^2}\, d \, \varepsilon_{n} 
\mathop{=}^{def} \delta \varepsilon_{n}  \quad \text{for all } n \geq 1.   
\end{displaymath}
We thus  get the desired conclusion \eqref{Main}  of Theorem~\ref{Principal}. \qed


 \section{An upper bound} \label{section 5}
 
We do not obtain a fairly good upper bound, and we shall content ourselves with the following result, whose proof is quite simple and, for the case  $\alpha = -1$, 
partly contained in \cite{Parf}, but under a very cryptic form which is not easy to decipher. 
 
 \begin{theoreme}\label{Ternary} 
Let $\varphi$ be a Schur function and $\alpha\geq -1$.Then,  we have for the approximation numbers of  
$C_\varphi \colon \mathfrak{B}_\alpha\to\mathfrak{B}_\alpha$  the upper bound:
\begin{equation}\label{Super} 
a_{n} (C_\varphi) \leq C \inf_{0 < h < 1} \bigg[n^{\frac{\alpha+1}{2}} (1 - h)^n 
+ \sqrt{\frac{\rho_{\varphi,\alpha+2} (h)}{h^{2+\alpha}}}\bigg], \qquad n = 1,2,\ldots 
\end{equation} 
where $C$ is a constant. In particular, if  $\frac{\rho_{\varphi,\alpha+2} (h)}{h^{2+\alpha}} \leq \e^{- h/ A(h)}$, where the function 
$A \colon [0,1] \to [0,1]$ is  increasing,  with $A (0) = 0$ and with inverse function $A^{-1}$, we have: 
\begin{equation}\label{Supper} 
\qquad \qquad \quad a_{n} (C_\varphi) \leq Cn^{\frac{\alpha+1}{2}} \e^{-n A^{- 1}(1/2n)}, \qquad n = 1, 2, \ldots.
\end{equation}
\end{theoreme}
\smallskip

The proof of  \eqref{Super} uses a contraction principle which was first proved for $\alpha=-1$ (\cite{MEMOIRS}) and $\alpha=0$ (\cite{TRANS}), but  is 
also valid  for any $\alpha\geq -1$, as follows from the forthcoming work \cite{PREPR}. 
\smallskip

To prove Theorem~\ref{Ternary}, it will be convenient to prove first the following simple lemma.
\begin{lemme} 
Let $n$ be a positive integer, $g\in \mathfrak{B}_\alpha$ and $f (z) = z^n g (z)$. Then, we have:
\begin{equation} \label{Division} 
\Vert g \Vert_\alpha \leq C n^{\frac{\alpha+1}{2}} \Vert f \Vert_\alpha.
\end{equation}
\end{lemme}
\noindent{\bf Proof.} Let $w_n =  \frac{n! \Gamma (2+\alpha)}{\Gamma (n + 2 + \alpha)}$ and $f (z) = \sum_{n=0}^\infty a_n  z^n$. We first observe that 
\begin{equation} \label{Quotient} 
\qquad \qquad \qquad \frac{w_k}{w_{k + n}} \leq C n^{\alpha+1}, \qquad \forall k \geq 0, \quad \forall n \geq 1.
\end{equation}
 Indeed, we have:
\begin{align*}
\frac{w_k}{w_{k + n}} 
& = \frac{k!} {(k+n)!} \frac{\Gamma (k + \alpha + 2 + n)} {\Gamma (k +\alpha + 2)} 
=\prod_{j=1}^n \frac{(k + j + \alpha + 1)} {(k + j)} 
\leq \prod_{j=1}^n \frac{j + \alpha + 1} {j} \\ 
& = \prod_{j=1}^n\bigg(1 + \frac{\alpha + 1} {j} \bigg) 
\leq \exp \bigg[ (\alpha + 1)\sum_{j=1}^n \frac{1}{j} \bigg] \leq C n^{\alpha+1}, 
\end{align*} 
 which proves \eqref{Quotient}. \par 
Now, if $f (z) = \sum_{k = n}^\infty a_k z^k$, we have $g (z) = \sum_{k=0}^\infty a_{k+n}  z^k$ so that, using \eqref{Quotient}:
\begin{displaymath}
\Vert g \Vert^2_\alpha = \sum_{k=0}^\infty \vert a_{k+n} \vert^2 w_k = \sum_{l=n}^\infty \vert a_{l} \vert^2 w_{l - n} 
\leq C n^{\alpha + 1} \sum_{l=n}^\infty \vert a_{l} \vert^2 w_l 
=  C n^{\alpha+1}\Vert f \Vert^2_\alpha , 
\end{displaymath}
proving \eqref{Division}.  \qed 
\par \medskip
 
We shall now majorize $a_{n+1} (C_\varphi)$, but provided that we change the constant $C$, this makes no difference with majorizing $a_{n} (C_\varphi)$. 
The choice of the approximating operator $R$ of rank $\leq n$ for $C_\varphi$ is quite primitive, but in counterpart we shall estimate 
$\Vert C_\varphi - R \Vert$ rather sharply.  We denote by $P_n$ the projection operator defined by $P_{n} f = \sum_{k=0}^{n - 1} \hat f (k) z^k$ and we 
take $R = C_\varphi \circ P_n$, i.e. if we have $f (z) =\sum_{k=0}^\infty \hat f (k) z^k \in \mathfrak{B}_\alpha$, then 
$R (f) =\sum_{k=0}^{n - 1} \hat f (k) \varphi^k$, so that $(C_\varphi - R) f = C_{\varphi} (r)$, with, making use of \eqref{Division}:
\begin{equation} \label{OK yet} 
r (z) = \sum_{k=n}^\infty \widehat f (k) z^k = z^{n} s (z), \quad \text{with} \quad 
\Vert s \Vert^2_\alpha \leq C n^{\alpha+1} \Vert r \Vert^2_\alpha ,  
\  \Vert r \Vert_\alpha \leq \Vert f \Vert_\alpha.
\end{equation}
Assume that $\Vert f \Vert_\alpha \leq 1$, fix $0 < h < 1$ and denote by $\mu_h$ the restriction of the measure $A_{\varphi,\alpha+2}$ to the annulus 
$1 - h < \vert z \vert \leq 1$. Then, we have:
\begin{align*}
\Vert (C_\varphi - R) f\Vert^2_\alpha  
& =\Vert C_{\varphi} (r) \Vert^2_\alpha = \int_{\overline{\D}} \vert r (z) \vert^2 \, dA_{\varphi,\alpha+2} (z) \\
& \leq (1 - h)^{2n} \int_{\vert z \vert \leq 1 - h} \vert s (z) \vert^{2} \, dA_{\varphi,\alpha+2} (z) \\
& \hskip 100pt + \int_{1 - h < \vert z \vert \leq 1} \vert r (z) \vert^2 \,dA_{\varphi,\alpha+2}(z) \\
& \leq (1 - h)^{2n} \int_{\D} \vert s (z) \vert^{2} \,dA_{\varphi,\alpha+2} (z) + \int_{\overline{\D}} \vert r (z) \vert^2 \, d\mu_h (z)  \\ 
& = (1 - h)^{2n} \Vert C_{\varphi} (s)\Vert^2_\alpha + \int_{\overline{\D}} \vert r (z) \vert^2 \,d\mu_h (z) \\
& \leq C \bigg[ (1 - h)^{2n} \Vert s \Vert^2_\alpha + \int_{\overline{\D}} \vert r (z) \vert^2 \,d\mu_h (z) \bigg] \\
& \leq C \bigg[ n^{{\alpha+1}} (1 - h)^{2n} + \sup_{0 < t \leq h} \frac{\rho_{\varphi,\alpha+2} (t)} {t^{2 + \alpha}} \bigg] 
\end{align*}
if we  use \eqref{OK yet}, as well as  \eqref{A priori} under the form 
\begin{displaymath}
\int_{\overline{\D}} \vert r (z) \vert^2 \,d\mu_h (z) 
\leq C \sup_{0 < t \leq h} \frac{\rho_{\varphi,\alpha+2} (t)} {t^{2 + \alpha}} \, \Vert r \Vert^2_\alpha, 
\end{displaymath}
and we know that $\Vert r \Vert_\alpha \leq \Vert f \Vert_\alpha \leq 1$.
\par\smallskip

To get rid of the supremum with respect to $t$, we make use of the following inequality, which holds for $h \leq 1 - \vert \varphi (0) \vert$ and 
$0 < \varepsilon \leq 1$:
\begin{equation}
\rho_{\varphi, \alpha + 2} (\varepsilon h) \leq C \varepsilon^{\alpha + 2} \rho_{\varphi, \alpha + 2} (h).
\end{equation}
For $\alpha = 0$ or $\alpha = - 1$, this follows respectively from \cite{MEMOIRS}, Theorem~4.19, p.~55, and from \cite{TRANS}, Theorem~3.1. The 
general case is proved in \cite{PREPR}. Setting $t = \varepsilon h$ for $0 < t \leq h$, this also reads  
$\frac{\rho_{\varphi,\alpha+2} (t)} {t^{\alpha+2}} \leq C \frac{\rho_{\varphi,\alpha+2} ( h)} {h^{\alpha+2}}$, and we can forget the supremum in $t$ in 
the previous inequalities. Taking square roots, we get the relation \eqref{Super}.\par
\smallskip

When $\frac{\rho_{\varphi,\alpha+2} (h)} { h^{2+\alpha}} \leq \e^{- h/ A(h)}$, let us take for $h$ the nearly optimal value  $h = A^{-1} (1/2n) $, so that 
$h/ A (h) = 2n h$. We then  have from \eqref{Super}, since $(1 - h)^{2n} \leq \e^{- 2n h}$:
\begin{displaymath}
a_{n + 1} (C_\varphi)^2 \leq \Vert C_\varphi - R\Vert^2_\alpha 
\leq C n^{\alpha + 1} [\e^{- 2n h} + \e^{- h/ A (h)} ] \leq 2C n^{\alpha+1} \e^{-2n A^{-1}(1/ 2n)},  
\end{displaymath}
proving \eqref{Supper}, and ending the proof of Theorem~\ref{Ternary}. \qed
\bigskip

Let us now indicate three corollaries, which improve results of \cite{JFA}, \cite{REVI}, and \cite{TRANS} respectively.
 
\begin{corollaire} \label{Imprecise} 
Suppose that $\rho_{\varphi,\alpha+2} (h) \leq C h^{(2+\alpha) \beta}$ for some $\beta > 1$. Then: 
\begin{displaymath}
a_{n} (C_\varphi) \leq C n^{- \frac{(\beta - 1) (\alpha + 2)} {2}} (\log n)^{\frac{(\beta - 1) (\alpha + 2)}{2}}.
\end{displaymath}
In particular, $C_\varphi$ belongs to the Schatten class $S_p = S_{p} (\mathfrak{B}_\alpha)$ for each 
$p >  \frac{2}{(\beta  - 1) (\alpha + 2)} \cdot$ 
\end{corollaire}
\noindent{\bf Proof.} Set $\gamma = (\beta - 1) (\alpha+2)/2$, $a = (\alpha+1)/2$, and $c = a + \gamma$. If we apply \eqref{Super} of 
Theorem~\ref{Ternary}  with the value $h = c \log n/ n$ which satisfies $n^{a} \e^{- n h} = n^{- \gamma}$, as well as the inequality 
$(1 - h)^n \leq \e^{- nh}$, we get:
\begin{displaymath}
a_{n} (C_\varphi) \leq C \bigg[ n^{- \gamma} + \Big( \frac{\log n}{n} \Big)^\gamma \bigg] \leq C \Big( \frac{\log n}{n} \Big)^\gamma , 
\end{displaymath}
ending the proof.\qed
\par\medskip

In \cite{JFA}, we had only the assertion on Schatten classes,  for the single value   $\alpha = - 1$, and not the upper bound for the individual approximation 
numbers $a_{n} (C_\varphi)$. \par
\smallskip

\begin{corollaire} \label{Onto} 
Let $(\varepsilon_n)$ a sequence of positive numbers which tends to $0$. Then, there exists a Schur function $\varphi$ with the following properties: \par
\smallskip
{\rm 1.} $\varphi \colon \D \to \D$ is surjective and $4$-valent;\par
\smallskip 
{\rm 2.} $a_{n} (C_\varphi) \leq C \e^{-n\varepsilon_n}$, $n = 1, 2, \ldots$
\par \smallskip 
\noindent 
In particular, we can get $a_{n} (C_\varphi) \leq C \e^{- \frac{n}{\log (n + 1)}}$ and $C_\varphi$ is in every Schatten class 
$S_{p} (\mathfrak{B}_\alpha)$, $p>0$.
\end{corollaire} 

Notice that the sequence $(\varepsilon_n)$ in the statement cannot be dispensed with. Indeed, if $\varphi$ is surjective, we surely have 
$\Vert \varphi \Vert_\infty = 1$! And we know from Theorem~\ref{Lastresult} that $\beta (C_\varphi) = 1$ in that case.\par \medskip

We begin with a lemma of independent interest.
\begin{lemme} \label{Arbitrary} 
Let $\delta \colon  (0,1] \to \R$ be a positive and non-decreasing function. Then there exists a Schur function $\varphi$ with the following properties: 
\par\smallskip
{\rm 1.} $\varphi \colon \D \to \D$ is surjective and $4$-valent; 
\par\smallskip 
{\rm 2.} $\rho_{\varphi, \alpha+2} (h) \leq \delta (h)$, for $h > 0$ small enough.
\end{lemme}

\noindent{\bf Proof.} We begin with the case $\alpha = - 1$. Set, for $a = 1/2$:
\begin{displaymath}
\Phi_{a} ( z) = \frac{a - z}{1 - az} \raise 1,5pt \hbox{,} \quad B = \Phi_{a}^2 ,
\end{displaymath}
and $C = \frac{1 + a}{2 (1 - a)} = 3/2$. Note that $B \big(\frac{2a}{a + 1} \big) = B (0)$. Let now 
\begin{displaymath}
b_n = \frac{1}{4n \pi} \raise 1,5pt \hbox{,} \qquad \varepsilon (h) = \frac{1}{2} \delta (h/C), \qquad  \varepsilon_n = \varepsilon (b_{n+1}) .   
\end{displaymath}
In  the proof of Theorem~4.1 of \cite{REVI}, using an argument of harmonic measure and of barrier, we have found a  $2$-valent  symbol $\varphi_1$ with 
$\varphi_1 (\D) =\D^*$ such that, noting $\rho_\varphi$ for $\rho_{\varphi,1}$:
\begin{equation} \label {Bound above} 
b_{n+1} < h \leq b_n \quad \Longrightarrow \quad \rho_{\varphi_{1}} (h) \leq \varepsilon_n.
\end{equation}
This gives $\rho_{\varphi_{1}} (h) \leq \varepsilon (b_{n+1}) \leq \varepsilon (h)$. Let now, as in \cite{REVI}, $\varphi = B \circ \varphi_1$. This Schur 
function is surjective (since $\varphi (\D) = B(\D^*) = B(\D) = \D$), and $4$-valent. Moreover, if $I = (u, v)$ is an arc of $\T$ of length $h < 1/ 2$ and 
$J = (\frac{u}{2}, \frac{v}{2})$, we have $B^{-1} (I) \subset \Phi_{a} (J) \cup \Phi_{a} (-J) = I_1\cup I_2$, where $I_1, I_2$ are two arcs of $\T$ of length 
at most $\Vert P_a\Vert_\infty (h/ 2) = C h$, since $\Phi_a$ being an inner function, we have (\cite{NORD}), $P_a$ being the Poisson kernel at $a$: 
\begin{displaymath}
m_{\Phi_a} = P_{a} m . 
\end{displaymath}
Hence, using \eqref{Bound above}, we obtain: 
\begin{displaymath}
m_{\varphi} (I) = m_{\varphi_1} (B^{-1} (I)) \leq m_{\varphi_1} (I_1) + m_{\varphi_1} (I_2) 
\leq 2 \rho_{\varphi_1} (Ch) \leq 2 \varepsilon (Ch) =\delta (h), 
\end{displaymath}
and $\rho_{\varphi} (h) \leq \delta (h)$ for small $h$, by passing to the supremum on all $I$'s.
\par\medskip

For the general case $\alpha\geq - 1$, we use the following extension of an inequality from \cite{TRANS} (which treats the case $\alpha = 0$, see 
Remark before Corollary~3.11):
\begin{lemme} 
For small $h$, namely $0 < h < (1 - \vert \varphi (0) \vert)/ 4$, we have, for every $\alpha > -1$:
\begin{equation} \label{Borrow} 
\rho_{\varphi, \alpha + 2} (h) \leq C [\rho_{\varphi} (Ch)]^{\alpha+2}.
\end{equation}
\end{lemme}
\noindent{\bf Proof.} Let us define, as in \cite{Shap}, the generalized Nevanlinna counting function $N_{\varphi,\alpha+2}$ by the formula
\begin{displaymath}
\qquad N_{\varphi, \alpha+2} (w) = \sum_{\varphi(z) = w} [\log (1/\vert z \vert )]^{\alpha+2}, 
\qquad w\in\D \setminus \{\varphi (0)\}. 
\end{displaymath}
The case $\alpha = -1$ corresponds to the usual Nevanlinna counting function, which will be denoted by $N_\varphi$. The partial Nevanlinna counting 
function $N_{\varphi} (r, w)$ is defined, for $0 \leq r \leq 1$, by:
\begin{displaymath}
N_{\varphi} (r, w) = \sum_{\varphi (z) = w} \log^+ (r/ \vert z\vert),
\end{displaymath}
so that $N_{\varphi} (1, w) = N_{\varphi} (w)$.\par
Since $\alpha + 2 \geq 1$, we have the obvious but useful inequality:
\begin{equation} \label{Obvious} 
N_{\varphi,\alpha+2} (w) \leq [N_{\varphi} (w)]^{\alpha+2}.
\end{equation}
We shall also make use of the following identity, due to J.~Shapiro (\cite{Shap}, Proposition~6.6, where a weight $1/ r$ is missing), and which can easily be 
checked after two integrations by parts:
\begin{equation}\label{Identity} 
N_{\varphi,\alpha+2} (w) = (\alpha + 2) (\alpha + 1) \int_0^1 N_{\varphi }(r, w) [\log (1/ r)]^{\alpha} \,\frac{dr}{r} \cdot 
\end{equation} 
As it was noticed in (\cite{TRANS}, Theorem~3.10), this formula reads, for $w$ close to the boundary, as follows, for 
$0 < h < (1 - \vert \varphi (0) \vert)/ 4$ and $\vert w \vert > 1 - h$:
\begin{equation}\label{Superidentity}
N_{\varphi,\alpha+2} (w) = (\alpha+2) (\alpha+1) \int_{1/3}^{1} N_{\varphi} (r, w) [\log (1/r)]^{\alpha} \, \frac{dr}{r} \cdot 
\end{equation} 
Under the same conditions on $h$ and $w$, this obviously implies:
\begin{displaymath}
N_{\varphi,\alpha+2} (w) \geq \frac{1}{C} \int_{1/3}^1 N_{\varphi} (r, w) (1 - r^2)^\alpha r\, dr 
= \frac{1}{C} \int_0^1 N_{\varphi} (r, w) (1 - r^2)^\alpha r\,dr.
\end{displaymath}
Now, using the same arguments as in  \cite{TRANS}, Theorem~3.10  and in particular using 
\eqref{Superidentity} for $\varphi_r (z) = \varphi (rz)$, the identity $N_{\varphi} (r, w) = N_{\varphi_r} (w)$ and an integration in polar coordinates, 
we get:  
\begin{equation}\label{Domination} 
\sup_{\vert w \vert\geq 1 - h} N_{\varphi,\alpha + 2} (w) \geq \frac{1}{C} \, \rho_{\varphi,\alpha+2} (h/C) .
\end{equation}
The end of the proof is easy: changing $h$ into $Ch$ and using successively \eqref{Domination} and \eqref{Obvious}, we get for small $h$, depending on 
$\varphi$:
\begin{displaymath}
\rho_{\varphi, \alpha + 2} (h) \leq C  \hskip -3pt \sup_{\vert w \vert \geq 1 - Ch} N_{\varphi, \alpha+2} (w) 
\leq C \hskip -3pt \sup_{\vert w \vert \geq 1 - Ch} [N_{\varphi} (w)]^{\alpha+2} 
\leq C\, [\rho_{\varphi} (Ch)]^{\alpha+2} , 
\end{displaymath}
the last inequality coming from \cite{MA}, Theorem~3.1. This ends the proof of \eqref{Borrow}. 
\par\bigskip

Going back to the proof of Lemma~\ref{Arbitrary}, if we apply the already settled case $\alpha = - 1$ to the function 
$\tilde \delta (h) = [\delta(h/ C)/C]^{\frac{1}{\alpha+2}}$, we obtain a surjective and $4$-valent Schur function $\varphi$ such that:
\begin{displaymath}
\rho_{\varphi, \alpha + 2} (h) \leq  C\, [\rho_{\varphi} (Ch)]^{\alpha+2}\leq  C\,[\tilde\delta (Ch)]^{\alpha+2} = \delta (h), 
\end{displaymath} 
for $h$ small enough.\qed
\par\bigskip

\noindent{\bf Proof of Corollary~\ref{Onto}.}  Set $a = (\alpha+1)/2$. Provided that we replace $(\varepsilon_n)$ by the decreasing sequence 
$(\varepsilon'_n)$ with $\varepsilon'_n = \frac{1}{n} + \sup_{k \geq n} \varepsilon_k \geq \varepsilon_n$, we can assume that $(\varepsilon_n)$ 
decreases. Let $A \colon [0,1]\to[0,1]$ be a function such that $A (0) = 0$, and which increases (as well as $A (t)/ t$) so slowly that 
$A (\varepsilon_n + a (\log n/ n)) \leq 1/ 2n$; therefore  $A^{-1} (1/ 2n) \geq \varepsilon_n + a (\log n/ n)$ and 
\begin{displaymath}
n^{a} \e^{- n A^{-1} (1/ 2n)} \leq \e^{- n \varepsilon_n} . 
\end{displaymath}
We now apply Lemma~\ref{Arbitrary} to the non-decreasing function $\delta (h) = h^{2+\alpha} \e^{- h/A (h)}$ to get the result, in view of  
\eqref{Supper} of Theorem~\ref{Ternary}. \qed 
\par\bigskip

Our last corollary involves Hardy-Orlicz spaces $H^\psi$ and Bergman-Orlicz spaces $\mathfrak{B}^\psi$. For the definitions, we refer to \cite{MEMOIRS}. 

\begin{corollaire}\label{Less} 
There exists a Schur function $\varphi$ and an Orlicz function $\psi$ such that $C_\varphi \colon H^\psi \to H^\psi$ is compact whereas 
$C_\varphi\colon \mathfrak{B}^\psi \to \mathfrak{B}^\psi$ is not compact. Moreover, the approximation numbers $a_{n} (C_\varphi)$ 
of $C_\varphi \colon \mathfrak{B}_\alpha \to \mathfrak{B}_\alpha$ satisfy the upper estimate $a_{n}(C_\varphi)\leq a \, \e^{- b \sqrt n}$ where 
$a$, $b$ are positive constants independent of $n$, and therefore $C_\varphi$  belongs to $\bigcap_{p > 0} S_p (\mathfrak{B}_\alpha)$.  
\end{corollaire}

\noindent{\bf Proof.}  Let $\alpha \geq - 1$ be fixed. The Schur function constructed in the proof of Theorem~4.2 of \cite{TRANS} satisfies the two first 
assertions, as well as $\rho_{\varphi} (h)/ h \leq \e^{- d/ h}$ for some positive constant $d > 0$. We now apply \eqref{Borrow} to get for small $h$:
\begin{displaymath}
\frac{\rho_{\varphi, \alpha+2} (h)} {h^{\alpha+2}} \leq C \, \frac{[\rho_{\varphi} (Ch)]^{\alpha+2}} {h^{\alpha+2}} 
\leq C^{\alpha+3} \e^{- (\alpha+2) d/ Ch} \leq a\, \e^{- b/ h}  
\end{displaymath}
for positive constants $a$ and $b$. We can thus apply  \eqref{Supper} of Theorem~\ref{Ternary}, for some $\delta > 0$, with the increasing function 
$A (h) = h^2/ \delta$ (hence $A^{-1} (x) = \sqrt{\delta x}$) to get the result, diminishing slightly $b$ to absorb the power factor 
$n^{\frac{\alpha+1}{2}}$. \qed 
\par\bigskip\goodbreak

\noindent{\bf Remark.} Let us alternatively consider the entropy numbers $e_n (C_\varphi)$ (see \cite{Carl} or \cite{K\"onig-livre}, page~69 for the definition) 
of composition operators. Those numbers are also a very good indicator of the ``degree of compactness'' of general operators $T \colon X\to Y$ where $X,Y$ are  
Banach spaces and are smaller than the approximation numbers, in the following weak sense (\cite{Pisier-livre}, page~64).
\begin{align}
& \qquad \sup_{1 \leq k \leq n} [k^{\alpha} e_{k} (T)] \leq C_\alpha  \sup_{1 \leq k \leq n} [k^{\alpha} a_{k} (T)], \qquad \forall \alpha > 0 .
\label{Banach weak} \\
& \qquad \quad (a_{n} (T)) \in \ell_q \quad \Longrightarrow \quad (e_{n} (T))\in\ell_q, \qquad  \quad \forall q > 0 .  \label{Banach strong}
\end{align}

The converse of  \eqref{Banach strong} does not hold in Banach spaces, but it does for operators between Hilbert spaces, by polar decomposition. More precisely, 
we have (\cite{Pisier-livre}, page~68) $a_{n} (T) \leq 4 \, e_{n} (T)$ and, in particular, $(e_{n} (T)) \in \ell_q$ if and only if $(a_{n} (T)) \in \ell_q$. \par
\medskip

We now have the following improved version of Theorem~\ref{Secondary}. Recall that 
$\varphi^{\#} (z) = \frac{\vert \varphi' (z) \vert (1 - \vert z \vert^2)} {1 - \vert \varphi (z) \vert^2}$ and $[\varphi] = \Vert \varphi^{\#} \Vert_\infty$.
\begin{theoreme} \label{Lastminute} 
Let $T = C_\varphi$ be a compact composition operator on $\mathfrak{B}_\alpha$, and $\gamma (T) = \liminf_{n\to\infty} [e_{n} (T)]^{1/n}$. Then:
\begin{equation} \label{Carl}
\gamma (T) \geq [\varphi]^{1/2}.
\end{equation}
\end{theoreme}
\noindent{\bf Proof.} We proceed as in the proof of Theorem~\ref{Secondary}. First, recall that the entropy numbers $e_{n} (T)$ also have the ideal 
property (\cite{K\"onig-livre}, page~69), namely:
\begin{displaymath} 
e_{n} (ATB) \leq \Vert A \Vert \, e_{n} (T) \,\Vert B \Vert. 
\end{displaymath} 
 Then, we use an improved Weyl-type inequality for entropy numbers, due to Carl and Triebel (\cite{CATR}), in which $(\lambda_{n} (T))_{n \geq 1}$ denotes 
the sequence of eigenvalues of $T$ rearranged in non-increasing order of  moduli and $C = \sqrt 2$:
 \begin{equation}\label{Improved} 
\Big(\prod_{k=1}^{n} \vert \lambda_{k} (T) \vert \Big)^{1/n} \leq C e_{n} (T).
\end{equation}
 It should be noted that this inequality can itself be improved (\cite{GOKOSC}): 
\begin{equation}\label{Reimproved} 
\Big(\prod_{k=1}^{n} a_{k} (T)\Big)^{1/n} \leq C e_{n} (T).
\end{equation}
Yet, the tempting similar inequality  $\Big(\prod_{k=1}^{n} \vert \lambda_{k} (T) \vert \Big)^{1/n} \leq C a_{n} (T)$ is wrong (even the inequality 
$\vert \lambda_{n} (T) \vert \leq C a_{n} (T)$ is wrong) as follows from an example of  (\cite{K\"onig-livre}, pages~133--134).  Note that \eqref{Reimproved} 
implies the following:
\begin{displaymath} 
\quad a_{n} (T) \geq \delta r^n \quad \Longrightarrow \quad e_{n} (T) \geq \frac{\delta}{C} r^{1/2} \, r^{n/2}.
\end{displaymath} 
This might explain why a square root appears in \eqref{Carl}, and tends to indicate that $[\varphi]$ should appear instead of $[\varphi]^2$ in 
Theorem~\ref{Secondary}. \par
\smallskip

Now, for  every $a\in\D$, let again $\Phi_a$ be defined by $\Phi_{a} (z) = \frac{a - z}{1 - \overline{a} z}$ \raise 1pt \hbox{,} for $z \in\D$. 
Set $b = \varphi (a)$ and define  $\psi = \Phi_{b} \circ \varphi \circ \Phi_a$.  We already know that $0$ is a fixed point of $\psi$  with  derivative  
$\psi' (0) = \phi^{\#}(a)$ and that $C_\psi = C_{\Phi_a} \circ C_\phi \circ C_{\Phi_{b}}$. We may assume that $\psi' (0) = \phi^{\#}(a) \neq 0$. The 
sequence of  eigenvalues of $C_\psi$  is then, as we have seen, $( (\psi' (0)^n )_{n \geq 0}$ (\cite{Shap-livre}, p.~96). The equation \eqref{Improved}  
then gives us, setting $r = \vert \psi' (0) \vert = \phi^{\#}(a)$: 
\begin{displaymath} 
e_{n} (C_\psi) \geq \frac{1}{C} \bigg(\prod_{k=0}^{n - 1} r^k \bigg)^{1/n} = \frac{1}{C} \, r^{(n - 1)/2}.
\end{displaymath} 
This clearly gives us $\gamma (C_\psi) \geq \sqrt r$. Now, since $C_{\Phi_a}$ and $ C_{\Phi_{b}}$ are invertible operators, the relation 
 $C_\psi = C_{\Phi_a} \circ C_\phi \circ C_{\Phi_{b}}$ and the ideal property of the numbers $e_n (T)$  imply that  
$\gamma (C_\varphi) = \gamma (C_\psi)$, and therefore, with the notation of \eqref{Dieze}, 
$\gamma(C_\varphi)\geq \big(\phi^{\#}(a)\big)^{1/2}$, for all $a\in\D$. Passing to the supremum on $a\in\D$,  we end the proof of  
Theorem~\ref{Lastminute}. \qed
 

\section{The explicit example of lens maps} \label{section 6}

To ease notation, we shall suppose in this section that $\alpha = - 1$, i.e. we are concerned with the Hardy space $H^2$. Fix $0 < \theta < 1$. Denote by 
$\mathbb H =\{z \in \C \,; \ \Re z > 0 \}$ the right half-plane, by $T \colon \D \to \C \setminus\{-1\}$ the involutive transformation defined by 
$T (z) = \frac{1 - z}{1 + z}$, which maps $\D$ to $\mathbb H$, and by $\tau_\theta$ the transformation 
$z \in {\mathbb H} \mapsto z^\theta \in {\mathbb H}$. Recall that the associated lens map $\varphi_\theta \colon \D \to \D$ is:
\begin{displaymath} 
\varphi_\theta = T\circ \tau_{\theta} \circ T. 
\end{displaymath} 
It is known that the associated composition operator on $H^2$ is in all Schatten classes $S_p$ (\cite{SHTA}, Theorem~6.3). Alternatively, one could use  
Luecking's criterion (\cite{LUECKING}). Therefore, its approximation numbers decrease rather quickly. Still more precisely, adapting techniques of 
Parfenov (\cite{Parf}, page~511), we might show the following (where $\beta_\theta, \gamma_\theta, \ldots$ are positive constants):
\begin{equation} \label{Nodetail} 
a_{n} (C_{\varphi_\theta}) \leq \gamma_\theta e^{-\beta_{\theta} \sqrt n}. 
\end{equation} 
We shall not detail this adaptation of Parfenov's methods from Carleson embeddings to composition operators, but shall dwell on the converse inequality, 
which is not proved in \cite{Parf}. First, the proof of the second assertion being postponed, we show that  there is no converse to the inequality of 
Theorem~\ref{Secondary}.

\begin{proposition}\label{Noconverse} 
The value of $[\varphi_\theta]$ for the lens map is 
\begin{equation}\label{First} 
[\varphi_\theta] = \theta.
\end{equation}
In particular,  $[\varphi_\theta]$ can be as small as we wish, although $\beta (C_{\varphi_\theta}) = 1$.
\end{proposition}

Recall that $\beta$ is defined in \eqref{beta} and $[\varphi]$ in \eqref{crochet}.\par
\medskip

\noindent{\bf Proof.}  First  note the simple   
\begin{lemme} 
Let $z \in\D$ and $v = T (z) \in {\mathbb H}$. Then:
\begin{displaymath} 
\vert T' (z) \vert (1 - \vert z \vert^2) = 2 \Re (T (z)) \quad \text{and} \quad \frac{\vert T' (v)\vert}{1 - \vert T (v) \vert^2} = \frac{1}{2 \Re v} \cdot
\end{displaymath} 
\end{lemme} 
The two equalities are the same because $\vert T' (v) \vert = \frac{1}{\vert T' (z) \vert}$ in view of $T = T^{-1}$. For the first one, we have:
\begin{displaymath} 
\vert T' (z) \vert ( 1 - \vert z \vert^2) = \frac{2 (1 - \vert z \vert^2)} {\vert 1 + z \vert^2} = 2 \Re (T (z)).
\end{displaymath} 
Let now $z \in \D$ and $w = T(z) \in {\mathbb H}$. By the chain rule, we have: 
\begin{displaymath} 
\varphi'_{\theta} (z) = T' (\tau_{\theta} (w)) \,\tau'_{\theta} (w) \,T'(z).
\end{displaymath} 
Taking moduli and using the lemma with $z$ and $v = \tau_{\theta} (w)$, we obtain:
\begin{displaymath} 
\frac{\vert \varphi'_{\theta} (z) \vert (1 - \vert z \vert^2)} {1 - \vert \varphi_{\theta} (z) \vert^2} 
= \frac{\vert T' (\tau_{\theta} (w)) \vert} {1 - \vert T (\tau_{\theta} (w)) \vert^2} \vert \tau'_{\theta} (w) \vert \, \vert T' (z) \vert (1 -\vert z \vert^2) 
= \frac{\vert\tau'_{\theta} (w) \vert \Re w} {\Re(\tau_{\theta} (w))} \cdot
\end{displaymath} 
Now,  setting $w = r \e^{it}$ with $r > 0$ and $- \pi/ 2 < t < \pi/2$, this writes as well:
\begin{displaymath} 
\varphi_{\theta}^{\#} (z) = \frac{\theta r^{\theta - 1} r \cos t} {r^{\theta} \cos \theta t} = \frac{\theta \cos t}{\cos \theta t}.
\end{displaymath} 
Using the fact that $w$ runs over $\mathbb H$ as $z$ runs over $\D$ and that the cosine decreases on $(0,\pi/2)$, we obtain \eqref{First} by taking $t=0$.
\qed
\par\medskip

We now prove the second assertion of Proposition~\ref{Noconverse} under the following form (the small roman and Greek letters 
$a_\theta, \ldots, \beta_\theta, \ldots$ will denote positive constants depending only on $\theta$):

\begin{proposition}\label{Ugly} 
There exist constants $b_\theta, c_\theta, , \beta_\theta, \gamma_\theta$ with \smash{$b_\theta = \pi \sqrt{\frac{2 (1 - \theta)}{\theta}}$} such that: 
\begin{equation}\label{Insane}
c_\theta \, \e^{-b_\theta \sqrt{n}} \leq a_{n} (C_{\varphi_\theta}) \leq \gamma_\theta \,\e^{- \beta_\theta \sqrt{n}}.
\end{equation}
In particular, we have  $\beta (C_{\varphi_\theta}) = 1$ and $C_{\varphi_\theta}$ is in all Schatten classes $S_p$, $p>0$ but its approximation numbers do not 
decrease exponentially.
\end{proposition} 

The upper bound is \eqref{Nodetail}. For the lower bound, we shall need  two simple lemmas.
\begin{lemme}\label{Hadlac} 
Let $0 < \sigma < 1$ and $u = (u_j)$ be a sequence of points of $\D$ such that $\frac{1 -\vert u_{j+1}\vert} {1 - \vert u_j \vert} \leq \sigma$. Then, the 
Carleson constant $\delta_u$ of the sequence $u$ satisfies: 
\begin{displaymath} 
\qquad \qquad \qquad \delta_u \geq \exp \bigg(- \frac{a}{1 - \sigma}\bigg), \qquad \text{with } a = \frac{\pi^2}{2} \cdot
\end{displaymath} 
\end{lemme}

\noindent{\bf Proof.} We use the following fact (\cite{Hoffman-livre}, pages~203--204):
\begin{equation}\label{Hoffman} 
\delta_u \geq \prod_{j=1}^\infty \bigg( \frac{1 - \sigma^j}{1 + \sigma^j} \bigg)^2 . 
\end{equation} 
This implies  $\log \delta_u \geq 2 \sum_{j=1}^\infty \log (\frac{1 - \sigma^j}{1 + \sigma^j})$. Now, expanding the logarithm in power series and permuting 
sums, we note that:
\begin{align*}
2 \sum_{j=1}^\infty  \log \Big(\frac{1 + \sigma^j}{1 - \sigma^j}\Big) 
& = 4 \sum_{k=0}^\infty \frac{\sigma^{2k+1}} {(2k + 1) (1 - \sigma^{2k+1})} \\  
& \leq 4 \sum_{k=0}^\infty \frac{1}{(2k+1)^2 (1 - \sigma)} = \frac{a}{1 - \sigma}  \raise 1,5pt \hbox{,}
\end{align*}
where we used 
$1 - \sigma^{2k+1} \geq (2k+1)(1 - \sigma) \sigma^{2k+1}$ and $\sum_{k=0}^\infty \frac{1}{(2k+1)^2} = \pi^2/8$. So that   
$\delta_u  \geq \exp \big(- a/ (1 - \sigma) \big)$, which was to be proved. \qed
\medskip

 The second lemma is similar. 
\begin{lemme}\label{Boue}  
Let $0 < \sigma < 1$,  $ u_j = 1 - \sigma^j$,  $v_j = \varphi_\theta (u_j)$ and  $v = (v_j)$. Then, the Carleson constant $\delta_v$ of the sequence $v$ satisfies:
\begin{displaymath} 
\qquad \qquad \qquad \delta_v \geq  \exp\Big(- \frac{a_\theta}{1 - \sigma}\Big) , \qquad \text{with } a_\theta = \frac{\pi^2}{2^{\theta} \theta} \cdot
\end{displaymath} 
\end{lemme}
\noindent{\bf Proof.}  We first note that 
$1 - \varphi_\theta (r) = \frac{2 (1 - r)^\theta}{(1 + r)^\theta + (1 - r)^\theta}$\raise 1pt \hbox{,}  and so  
\begin{displaymath} 
\frac{1 - v_{j+1}} {1  -v_j} 
= \sigma^\theta \frac{\sigma^{j \theta} + (2 - \sigma^j)^{\theta}} {\sigma^{(j+1)\theta} + (2 - \sigma^{j+1})^\theta} = \sigma_j , 
\end{displaymath} 
with $\sigma_j \leq \sigma' = 1 - \frac{\theta}{2} 2^\theta (1 - \sigma)$. To see this, observe that: 
\begin{displaymath} 
1 - \sigma_j = \frac{(2 - \sigma^{j+1})^\theta - (2 \sigma - \sigma^{j+1})^\theta} { \sigma^{(j+1)\theta} + (2 - \sigma^{j+1})^{\theta}}
\mathop{=}^{def} \frac{N}{D} \geq \theta 2^{\theta - 1} (1 - \sigma) = 1 - \sigma' . 
\end{displaymath} 
Indeed, the function $f (x) = x^\theta + (2 - x)^\theta$ increases on $[0, 1]$, so $D \leq f (1) = 2$. On the other hand, the mean-value theorem gives 
$N = 2 (1 - \sigma) \theta \, c^{\theta - 1} \geq \theta (1 - \sigma) 2^\theta$ for some $c \in (0, 2)$. Lemma~\ref{Hadlac} then gives the result for the 
sequence $v$. \qed
\medskip

\noindent{\bf Proof of  Proposition~\ref{Ugly}.} Fix an integer $n \geq 1$, and take $(u_j)$, $(v_j)$ as in Lemma~\ref{Boue}. We have 
$\varphi_\theta (0) = 0$, $\vert \varphi_\theta (z) \vert \leq \vert z \vert$ and so for $0 < r <1$:
\begin{displaymath} 
\frac{1 - r^2}{1 - \varphi_\theta (r)^2} \geq \frac{1 -  r}{1 - \varphi_\theta (r)} = 
\frac{(1 - r)[(1 - r)^\theta + (1 + r)^\theta]} {2 (1 - r)^\theta} \geq \frac{(1 - r)^{1 - \theta}}{2} \raise 1,5 pt \hbox{,} 
\end{displaymath} 
implying 
\begin{displaymath} 
\qquad \qquad \frac{1 - u_{j}^2} {1 - v_{j}^2} \geq \frac{1}{2} \, \sigma^{n (1 - \theta)} , \quad \text{for } 1 \leq j \leq n.
\end{displaymath} 

Let now $R$ be an operator of rank $< n$. There exists a function $f = \sum_{j=1}^n \lambda_j K_{u_j} \in H^2 \cap \ker R$ with $\Vert f \Vert =1$. 
We thus have, denoting by $C_u$ and $C_v$ the interpolation constants of the sequences $u$ and $v$, and using Lemma~\ref{Riesz basis} twice: 
\begin{align*}
\Vert C_{\varphi_\theta}^* - R \Vert^2 
& \geq \Vert C_{\varphi_\theta}^* (f) - R (f) \Vert^2 
= \Vert C_{\varphi_\theta}^* (f) \Vert^2 = \bigg\Vert \sum_{j=1}^n \lambda_j K_{v_j} \bigg\Vert^2 \\
& \geq C_{v}^{- 2} \sum_{j=1}^n \vert \lambda_{j} \vert^2 \Vert K_{v_j} \Vert^2
= C_{v}^{- 2} \sum_{j=1}^n \frac{\vert \lambda_{j} \vert^2}{1 -  v_j^2} \\
& \geq \frac{1}{2} \, C_{v}^{-2} \sigma^{n (1 - \theta)} \sum_{j=1}^n \frac{\vert \lambda_{j} \vert^2}{1 -  u_j^2}
\geq \frac{1}{2} \,C_{u}^{-2} C_{v}^{-2} \sigma^{n (1 - \theta)} \Vert f \Vert^2 \\
& = \frac{1}{2}\, C_{u}^{-2} C_{v}^{-2} \sigma^{n (1 - \theta)}.
\end{align*}
Therefore, $a_{n} (C_{\varphi_\theta}) = a_{n} (C_{\varphi_\theta}^*) \geq \frac{1}{2}\,C_{u}^{-1} C_{v}^{-1} \sigma^{n (1 - \theta)/2}$.  
But it follows from \eqref{Interpolation constant}, Lemma~\ref{Hadlac} and Lemma~\ref{Boue} that  $C_u$, $C_v$ satisfy, provided that we now take the 
value $a_\theta = \frac{\pi^2}{\theta} > \frac{\pi^2}{2} + \frac{\pi^2}{2^{\theta} \theta}$, since $\theta + 2^{1 - \theta} < 2$, to absorb the logarithmic 
factor of \eqref{Interpolation constant}:
\begin{displaymath} 
C_u C_v \leq c_{\theta}^{-1} \exp\big(a_\theta / (1 - \sigma)\big). 
\end{displaymath} 
The preceding now gives us ($c_\theta$ changing from line to line):
\begin{displaymath} 
a_{n} (C_{\varphi_\theta}) \geq c_\theta \exp\Big(- \frac{a_\theta}{1 - \sigma}\Big) \exp\bigg(\frac{n (1 - \theta)}{2} \log \sigma\bigg).
\end{displaymath} 
Finally, adjust $\sigma = 1 - \lambda n^{- 1/2}$ so that $\frac{a_\theta}{\lambda} = \frac{1 - \theta}{2} \lambda$, i.e. 
$\lambda = \sqrt{\frac{2a_\theta}{1 - \theta}}$ and use $\log (1 - x) \geq - x - x^2$ for $0 \leq x \leq 1/2$; this gives \eqref{Insane} with the value 
\begin{displaymath} 
b_\theta = \frac{2 a_\theta}{\lambda} = \sqrt{2 a_{\theta}(1 - \theta)} = \pi \sqrt{\frac{2(1 - \theta)}{\theta}} \raise 1,5pt \hbox{,}
\end{displaymath} 
and that ends the proof of Proposition~\ref{Ugly}. \qed
\bigskip

\noindent{\bf Remarks.} \par
1) The procedure used here to get lower estimates for the approximation numbers for lens maps might be easily adapted to a general symbol, to provide a new 
proof of Theorem~\ref{Secondary}. But the value of $\beta (C_\varphi)$ which we obtain in the general case is worse than the one obtained in 
Section~\ref{section 3}, therefore we did not think it useful to include this second proof. \par

2) It is easy to see that, for the lens map $\varphi_\theta$, one has $\rho_{\varphi_\theta} (h) \approx h^{1/\theta}$. Then 
Corollary~\ref{Imprecise} gives $a_{n} (C_{\varphi_\theta}) \leq C \, n^{- \frac{1 - \theta}{2\theta}} (\log n)^{\frac{1 - \theta}{2\theta}}$ and so 
$C_{\varphi_\theta}\in S_p$ for all $p > 2 \theta/ (1 - \theta)$. On the other hand, we know (\cite{SHTA}) that 
$C_{\varphi_\theta} \in \bigcap_{p > 0} S_p$, so that $a_{n} (C_{\varphi_\theta})$ must be rapidly decreasing: 
$a_{n} (C_{\varphi_\theta}) \leq C_{q} n^{-q}$ for all $q>0$. This shows that Theorem~\ref{Ternary} is very imprecise in general, becoming more accurate 
when $\rho_\varphi$ is very small, as this is the case in Corollary~\ref{Less}. \par \smallskip 

We hope to return to upper bounds for  approximation numbers in another work. 

\bigskip

\noindent{\bf Acknowledgement.} Part of this work was made during a visit of the second named author at the Departamento de An\'alisis Matem\'atico of 
the Universidad of Sevilla in April 2011; it is a pleasure to thank all people of this department for their warm hospitality. The third named author is partially 
supported by a Spanish research project MTM 2009-08934.


\vbox{\noindent {\small \it  {\rm Daniel Li}, Univ Lille Nord de France, 
U-Artois, Laboratoire de Math\'ematiques de Lens EA~2462, 
F\'ed\'eration CNRS Nord-Pas-de-Calais FR~2956, 
Facult\'e des Sciences Jean Perrin, 
Rue Jean Souvraz, S.P.\kern 1mm 18, 
F-62\kern 1mm 300 LENS, FRANCE, 
daniel.li@euler.univ-artois.fr
\\
{\rm Herv\'e Queff\'elec}, Univ Lille Nord de France, 
USTL, Laboratoire Paul Painlev\'e U.M.R. CNRS 8524, 
F-59\kern 1mm 655 VILLENEUVE D'ASCQ Cedex, 
FRANCE,  
queff@math.univ-lille1.fr
\\
{\rm Luis Rodr{\'\i}guez-Piazza}, Universidad de Sevilla, 
Facultad de Matem\'aticas, Departamento de An\'alisis Matem\'atico,  
Apartado de Correos 1160, 
41\kern 1mm 080 SEVILLA, SPAIN, 
piazza@us.es}
}

\end{document}